\def\doneprime{\end{document}}
\expandafter\ifx\csname bdtMiddleFlag\endcsname\relax
  \let\done=\doneprime
  \documentclass{article}
\catcode`\@=11  
\def\@startline{\global\@curtabmar\@nxttabmar\relax
   \global\@curtab\@curtabmar\setbox\@curline\hbox
    {}\@startfield\global\lifalse\strut}
\catcode`\@=12  

\newcommand{\For}{{\bf for} }

\newcommand{\If}{{\bf if} }

\newcommand{\Return}{{\bf return} }
\newcommand{\Comment}{$\!\triangleright$ }

\newcommand{\proc}[1]{
\ifmmode\mathord{\mathcode`-="702D\sc#1\mathcode`\-="2200}\else{\sc#1}\fi}

\newcommand{\func}[1]{\ifmmode\mathop{\mathcode`-="702D\rm #1\mathcode`\-="2200}\nolimits\else$\mathop{\mathcode`-="702D\rm #1\mathcode`\-="2200}\nolimits$\fi}

\newif\ifcodeinbox

\newcounter{codelinenumber}
\newcommand{\zeroli}{\setcounter{codelinenumber}0}

\catcode`\@=11  
\def\@startline{\global\@curtabmar\@nxttabmar\relax
   \global\@curtab\@curtabmar\setbox\@curline\hbox
    {}\@startfield\global\lifalse\strut}
\catcode`\@=12  

\newenvironment{code}{\global\codeinboxtrue
\setbox\strutbox\hbox{\vrule height 9pt depth 4pt width0pt}
\noindent\begin{tabbing}
\zeroli\setlength{\tabbingsep}{1em}
\hspace*{1em}\=999\ifdoubledigit9\fi
\=\ {\bf if} \={\bf then} \={\bf if} \={\bf then}
	\={\bf if} \={\bf then} \={\bf if} \={\bf then} \={\bf if} \={\bf then}
	\={\bf if} \=\+\+\kill}{\end{tabbing}\global\codeinboxfalse}

\catcode`\@=11  
\newif\ifdoubledigit
\newcommand{\codebox}[1]{\setbox0=\vbox{\begin{code}#1\end{code}}
\ifnum\c@codelinenumber>9\global\doubledigittrue\else\doubledigitfalse\fi
\vskip1sp\noindent\hskip-14pt
\parbox{\textwidth}{\begin{code}\protect#1\end{code}}\global\let\@currentlabel=\thechapter}
\catcode`\@=12  

\catcode`\@=11  
\newif\ifli
\newcommand{\li}{\global\litrue\stepcounter{codelinenumber}
\ifdoubledigit
\hbox to8pt{\hss\thecodelinenumber\hskip5pt}\else
\hbox to8pt{\hskip-1pt\thecodelinenumber\hss}\fi
\xdef\@currentlabel{\p@codelinenumber\thecodelinenumber}\'}
\catcode`\@=12  

\makeatletter
\def\fnum@figure{{\bf Figure \thefigure}}
\def\fnum@table{{\bf Table \thetable}}
\let\@mycaption\caption

\newcounter{totalcaptions}
\newcounter{totalart}

\newcommand{\psfigures}{0}

\newcommand{\macfigures}{\renewcommand{\psfigures}{2}}

\macfigures	

\newlength{\halffigspace} \newlength{\wholefigspace}
\newlength{\figruleheight} \newlength{\figgap}
\newcommand{\setfiglengths}{\ifnum\psfigures=1\setlength{\figruleheight}{\hruleheight}\setlength{\figgap}{1em}\else\setlength{\figruleheight}{0pt}\setlength{\figgap}{0em}\fi}
\newcommand{\figspace}[2]{\ifnum\psfigures=0\leavefigspace{#1}\else
\setfiglengths
\setlength{\wholefigspace}{#1}\setlength{\halffigspace}{.5\wholefigspace}
\rule[-\halffigspace]{\figruleheight}{\wholefigspace}\hspace{\figgap}#2\fi}
\newlength{\widefigspacewidth}
\newcommand{\widefigspace}[2]{
\ifnum\psfigures=0\leavefigspace{#1}\else
\setfiglengths
\setlength{\widefigspacewidth}{28pc}
\addtolength{\widefigspacewidth}{-\figruleheight}
\setlength{\wholefigspace}{#1}\setlength{\halffigspace}{.5\wholefigspace}
\makebox[\widefigspacewidth][r]{#2\hspace{\figgap}}\rule[-\halffigspace]{\figruleheight}{\wholefigspace}\fi}
\newcommand{\leavefigspace}[1]{\setlength{\wholefigspace}{#1}\setlength{\halffigspace}{.5\wholefigspace}\rule[-\halffigspace]{0em}{\wholefigspace}}

\newlength{\macfigfill}
\makeatother
\newlength{\bbx}
\newlength{\bby}
\newcommand{\macfigure}[5]{\addtocounter{totalart}{1}
\ifnum\psfigures=2
\setlength{\bbx}{#2}\addtolength{\bbx}{#4}
\setlength{\bby}{#3}\addtolength{\bby}{#5}
\begin{flushleft}
\ifdim#4>28pc\setlength{\macfigfill}{#4}\addtolength{\macfigfill}{-28pc}\hspace*{-\macfigfill}\fi
\mbox{\psfig{figure=./#1.ps,
bbllx=#2,bblly=#3,bburx=\bbx,bbury=\bby}}
\end{flushleft}
\else\ifdim#4>28pc\widefigspace{#5}{#1}\else\figspace{#5}{#1}\fi\fi}
\makeatletter

\newcommand{\defi}[1]{{\em #1}\glossary{#1}}
\newcommand{\indx}[1]{$#1$\index{#1}}

\catcode`\@=11  
\def\@begintheorem#1#2{\sl \trivlist \item[\hskip \labelsep{\bf #1\ #2}]}
\def\pem{\ifdim \fontdimen\@ne\font >\z@ \rm \else \sl \fi}
\catcode`\@=12  

\newtheorem{proposition}{Proposition}[section]
\newtheorem{theorem}[proposition]{Theorem}
\newtheorem{lemma}[proposition]{Lemma}
\newtheorem{corollary}[proposition]{Corollary}
\newtheorem{porism}[proposition]{Porism}

\newtheorem{definition}[proposition]{Definition}

\newtheorem{xample}{Example}[section]

\newenvironment{example}{\begin{xample}\rm}{\end{xample}}
\newcommand{\prf}{\par{\noindent\em Proof.}}
\newenvironment{proof}{\prf \setlength{\parindent}{3ex}}{\qed\par\noindent}
\def\squarebox#1{\hbox to #1{\hfill\vbox to #1{\vfill}}}
\newcommand{\qedbox}{\vbox{\hrule\hbox{\vrule\squarebox{.667em}\vrule}\hrule}}
\newcommand{\qed}{\nopagebreak\mbox{}\hfill\qedbox\smallskip}

\newcommand{\thought}[1]{XXXXXX #1 XXXXXX}

\newcommand{\noaside}[1]{}

\newcommand{\tensZ}{\mathop{\otimes_\integers}}
\newcommand{\tens}{\otimes}

\newcommand{\dUnion}{\biguplus}

\renewcommand{\k}{{\bf k}} 

\newcommand{\Q}{\rationals}
\newcommand{\Z}{\integers}
\newcommand{\N}{\naturals}

\newcommand{\integers}{{\bf Z}}
\newcommand{\naturals}{{\bf N}}

\newcommand{\rationals}{{\bf Q}}

\newcommand{\implies}{\Rightarrow}

\newcommand{\C}[1]{{\protect\cal #1}}
\newcommand{\mbf}[1]{{\protect\underline{#1}}}

\newcommand{\smush}[1]{{\makebox[0pt][l]{$#1$}}}

\newcommand{\pl}{pl}
\def\gl_#1{ {\pl}_{\C #1} }

\newcommand{\LP}[2]{\left[{\cal #1} \mid {\cal #2} \right]}

\newcommand{\Sym}{{{\cal S}ym}}

\newcommand{\Div}{{{\cal D}iv}}

\newcommand{\Wedge}{\Lambda}
\newcommand{\superrat}{ {{\cal S}uper_\rationals} }
\newcommand{\super}[2]{{{\cal S}uper}^{#1}(#2)}
\newcommand{\Super}[1]{{{\cal S}uper}(#1)}
\newcommand{\Superlp}[2]{\Super{\LP{#1}{#2}}}
\def\superlp#1#2^#3{\super{#3}{\LP{#1}{#2}}}

\newcommand{\SuperQ}[1]{{ \superrat}(#1)}
\newcommand{\SuperQlp}[2]{\SuperQ{\LP{#1}{#2}}}
\def\superQlp#1#2^#3{\superrat{#3}{\LP{#1}{#2}}}
\def\Superembed#1#2{\superrat
   (\{x_{a,d}\}_{a\in\C {#1},\, d\in\C {#2}} )} 
\def\c(#1)!{{\bf c}\left(#1\right)!}

\newcommand{\lspanZ}{\mathop{{\rm span}_\Z}}

\newcommand{\diag}{{\prec_{\hbox{\protect\tiny diag}}}}

\def\D_#1{{\rm D}_{#1}}

\def\auxilliaryR_#1,#2^#3_{\mathop{\,\mathstrut_{#1,#2}{\rm R}^{#3}}}
\def\R_#1{\auxilliaryR_#1^_}

\def\Rd_#1^#2{\auxilliaryR_#1^#2_}

\newcommand{\gp}{\mathrel{\raise1pt\hbox{$\scriptscriptstyle+$}\kern-2.5pt\hbox{$>$}}}
\newcommand{\gm}{\mathrel{\raise1pt\hbox{$\scriptscriptstyle-$}\kern-2.5pt\hbox{$>$}}}
\newcommand{\lp}{\mathrel{\hbox{$<$}\kern-2.5pt\raise1pt\hbox{$\scriptscriptstyle+$}}}
\newcommand{\lm}{\mathrel{\hbox{$<$}\kern-2.5pt\raise1pt\hbox{$\scriptscriptstyle-$}}}

\newcommand{\Deruyts}{{Der^-}}

 \def\tabmod_#1{{\C S}^{#1}}

\newcommand{\Tab}{{\hbox{Tab}}}

\newcommand{\inwt}[2]{{{\rm init}_{#1}(#2)}}

\def\init_#1(#2){{\rm init}_{#1}\left(#2\right)}
\def\tail_#1(#2){{\rm tail}_{#1}\left(#2\right)}
\def\LT_#1(#2){{\rm LT}_{#1}\left(#2\right)}
\def\initSet_#1(#2){{\rm LT}^{\rm monom}_{#1}\left(#2\right)}

\newcommand{\searrowminus}{\hbox{$\searrow$}\kern-3.5pt\raise1.5pt\hbox{\hbox{\vrule width2pt height0.2pt}}}
\newcommand{\nwarrowplus}{\hbox{$\nwarrow$}\kern-3.5pt\raise1pt\hbox{\raise1pt\hbox{\vrule width2pt height0.2pt}
\kern-1pt\hbox{\vrule width0.2pt height2pt}}}
\newcommand{\uparrowplus}
{\hbox{$\uparrow$}\kern-2pt\raise1pt
\hbox{\raise1pt\hbox{\vrule width2pt height0.2pt}
\kern-1pt\hbox{\vrule width0.2pt height2pt}}}
\newcommand{\uparrowminus}{\hbox{$\uparrow$}
\kern-2pt\raise2.0pt\hbox{\hbox{\vrule width2pt height0.2pt}}}
\newcommand{\downarrowplus}
{\hbox{$\downarrow$}\kern-2pt\raise1pt
\hbox{\raise1pt\hbox{\vrule width2pt height0.2pt}
\kern-1pt\hbox{\vrule width0.2pt height2pt}}}
\newcommand{\downarrowminus}{\hbox{$\downarrow$}
\kern-2pt\raise2.0pt\hbox{\hbox{\vrule width2pt height0.2pt}}}
\newcommand{\nearrowminus}{\hbox{$\nearrow$}\kern-6pt\raise1.5pt\hbox{\hbox{\vrule width2pt height0.2pt}}}
\newcommand{\swarrowplus}{\hbox{$\swarrow$}\kern-6pt\raise1pt\hbox{\raise1pt\hbox{\vrule width2pt height0.2pt}
\kern-1pt\hbox{\vrule width0.2pt height2pt}}}

\newcommand{\neam}{\smash{\smush{{\kern2pt\raise3pt\hbox{\nearrowminus}}}}}
\newcommand{\swap}{\smash{\smush{{\kern2pt\raise3pt\hbox{\swarrowplus}}}}}
\newcommand{\nwap}{\smash{\smush{{\kern-12pt\raise3pt\hbox{\nwarrowplus}}}}}
\newcommand{\seam}{\smash{\smush{{\kern-12pt\raise3pt\hbox{\searrowminus}}}}}
\newcommand{\uap}{\smash{\smush{\kern-2pt\raise3pt\hbox{\uparrowplus}}}}
\newcommand{\uam}{\smash{\smush{\kern-2pt\raise3pt\hbox{\uparrowminus}}}}
\newcommand{\dap}{\smash{\smush{\kern-2pt\raise3pt\hbox{\downarrowplus}}}}
\newcommand{\dam}{\smash{\smush{\kern-2pt\raise3pt\hbox{\downarrowminus}}}}
\newcommand{\ua}{\smash{\smush{\kern-2pt\raise3pt\hbox{$\uparrow$}}}}
\newcommand{\da}{\smash{\smush{\kern-2pt\raise3pt\hbox{$\downarrow$}}}}

\newcommand{\bsearrowminus}
{\begin{picture}(20,20)(0,10)
\thinlines
\put(17,6){\makebox[0pt]{$\scriptstyle -$}}
\put(0,20){\vector(1,-1){20}}
\end{picture}}

\newcommand{\bnearrowminus}
{\begin{picture}(20,20)(0,10)
\thinlines
\put(14,6){\makebox[0pt]{$\scriptstyle -$}}
\put(0,0){\vector(1,1){20}}
\end{picture}}
\newcommand{\bswarrowplus}
{\begin{picture}(20,20)(0,10)
\thinlines
\put(14,6){\makebox[0pt]{$\scriptstyle +$}}
\put(20,20){\vector(-1,-1){20}}
\end{picture}}
\newcommand{\bnwarrowplus}
{\begin{picture}(20,20)(0,10)
\thinlines
\put(18,6){\makebox[0pt]{$\scriptstyle +$}}
\put(20,0){\vector(-1,1){20}}
\end{picture}}
\newcommand{\bnwarrow}
{\begin{picture}(20,20)(0,10)
\thinlines
\put(18,6){\makebox[0pt]{$\scriptstyle $}}
\put(20,0){\vector(-1,1){20}}
\end{picture}}

\newcommand{\buparrowplus}
{\begin{picture}(20,20)(0,10)
\thinlines
\put(14,8){\makebox[0pt]{$\scriptstyle +$}}
\put(10,0){\vector(0,1){20}}
\end{picture}}
\newcommand{\buparrowminus}
{\begin{picture}(20,20)(0,10)
\thinlines
\put(14,8){\makebox[0pt]{$\scriptstyle -$}}
\put(10,0){\vector(0,1){20}}
\end{picture}}

\newcommand{\bdownarrowplus}
{\begin{picture}(20,20)(0,10)
\thinlines
\put(14,8){\makebox[0pt]{$\scriptstyle +$}}
\put(10,20){\vector(0,-1){20}}
\end{picture}}
\newcommand{\bdownarrowminus}
{\begin{picture}(20,20)(0,10)
\thinlines
\put(14,8){\makebox[0pt]{$\scriptstyle -$}}
\put(10,20){\vector(0,-1){20}}
\end{picture}}
\newcommand{\buparrow}
{\begin{picture}(20,20)(0,10)
\thinlines
\put(10,0){\vector(0,1){20}}
\end{picture}}
\newcommand{\bdownarrow}
{\begin{picture}(20,20)(0,10)
\thinlines
\put(10,20){\vector(0,-1){20}}
\end{picture}}

\newcommand{\bneam}{\smash{\smush{{\kern3pt\raisebox{3pt}{\makebox[0pt][l]{\bnearrowminus}}}}}}
\newcommand{\bswap}{\smash{\smush{{\kern3pt\raisebox{3pt}{\makebox[0pt][l]{\bswarrowplus}}}}}}
\newcommand{\bnwap}{\smash{\smush{{\kern-3pt\raisebox{3pt}{\makebox[0pt][r]{\bnwarrowplus}}}}}}
\newcommand{\bnwa}{\smash{\smush{{\kern-3pt\raisebox{3pt}{\makebox[0pt][r]{\bnwarrow}}}}}}
\newcommand{\bseam}{\smash{\smush{{\kern-3pt\raisebox{3pt}{\makebox[0pt][r]{\bsearrowminus}}}}}}
\newcommand{\buap}{\smash{\smush{\kern0pt\raisebox{3pt}{\makebox[0pt]{\buparrowplus}}}}}
\newcommand{\buam}{\smash{\smush{\kern-2pt\raisebox{3pt}{\makebox[0pt][c]{\buparrowminus}}}}}
\newcommand{\bdap}{\smash{\smush{\kern-2pt\raisebox{3pt}{\makebox[0pt][c]{\bdownarrowplus}}}}}
\newcommand{\bdam}{\smash{\smush{\kern-2pt\raisebox{3pt}{\makebox[0pt][c]{\bdownarrowminus}}}}}
\newcommand{\bua}{\smash{\smush{\kern-2pt\raisebox{3pt}{\makebox[0pt][c]{$\buparrow$}}}}}
\newcommand{\bda}{\smash{\smush{\kern-2pt\raisebox{3pt}{\makebox[0pt][c]{$\bdownarrow$}}}}}

\newcommand{\word}[1]{\smush{\cdots #1 \cdots}}
\newcommand{\aword}[1]{{\cdots #1 \cdots}}

\newcommand{\bu}{\bullet}
\newcommand{\lprn}{\smash{\kern-5pt\smush{\left(\rule{0pt}{25pt}\right.}}\hfill}
\newcommand{\rprn}{\hfill\smash{\smush{\left.\rule{0pt}{25pt}\kern-5pt\right)}}}
\newcommand{\up}[1]{{\raisebox{15pt}{\makebox[0pt]{$#1$}}}}
\newcommand{\down}[1]{{\raisebox{-15pt}{\makebox[0pt]{$#1$}}}}

\newcommand{\sst}{\scriptstyle}

\newlength{\across}
\newlength{\sqlength}
\newlength{\vstrutht}
\newlength{\vstrutdp}
\newlength{\vstrutwd}
\def\tableau#1#2#3{\begingroup
\renewcommand{\arraystretch}{0}
\renewcommand{\vbar}{\llap{\makebox[4pt][l]{\smash{\vstrut}}}}
\settowidth{\sqlength}{$#1$}
\setlength{\across}{#2\sqlength}
\newcommand{\struttext}{\strut} 
\def\abox##1{\makebox[\sqlength][l]{\vphantom{\struttext} $\,##1$}}
\def\closerow{\relax}
\def\boxify##1,##2 {\def\testarg{##1}\def\comparg{\closerow}
        \ifx\comparg\testarg\let\next=\relax\else\abox{##1}\phantom{##2}
        \let\next=\boxify\fi\next}
\def\sboxify##1,##2 {\def\testarg{##1}\def\comparg{\closerow}
\ifx\comparg\testarg\let\next=\relax\else\mbox{
\vphantom{\struttext} $\,##1$}\,\phantom{##2}
\let\next=\sboxify\fi\next}
\def\srow##1{\mbox{\sboxify ##1, \closerow, }}
\def\comparg{f}  
\def\compargc{c} 
\def\compargs{s} 
\def\compargt{t} 
\def\compargw{w} 
\def\testarg{#3}
\ifx\testarg\comparg
\addtolength{\across}{.2cm}   
\def\row##1{\framebox[\across][l]{\boxify ##1, \closerow, }}
\def\vrow##1##2##3{\makebox[##1\sqlength][l]{}
                     \framebox[##2\sqlength][l]{\boxify ##3, \closerow, }}
\def\irow##1##2##3{\makebox[##1\sqlength][l]{}
                     \makebox[##2\sqlength][l]{\boxify ##3, \closerow, }}
\setlength{\vstrutht}{8.5pt} \addtolength{\vstrutht}{1.5\fboxsep} 
\setlength{\vstrutdp}{3.5pt} \addtolength{\vstrutdp}{1.5\fboxsep}
\setlength{\vstrutwd}{1pt} 
\else\ifx\testarg\compargc
\addtolength{\across}{.2cm}   
\def\|{\hfill\vrule}
\def\row##1{\frame{\makebox[\across][l]{a\boxify ##1, \closerow, }}}
\def\vrow##1##2##3{\makebox[##1\sqlength][l]{}
                     \framebox[##2\sqlength][l]{\boxify ##3, \closerow, }}
\def\irow##1##2##3{\makebox[##1\sqlength][l]{}
                     \makebox[##2\sqlength][l]{\boxify ##3, \closerow, }}
\setlength{\vstrutht}{8.5pt} \addtolength{\vstrutht}{1.5\fboxsep} 
\setlength{\vstrutdp}{3.5pt} \addtolength{\vstrutdp}{1.5\fboxsep}
\setlength{\vstrutwd}{1pt} 
\else\ifx\testarg\compargs
\renewcommand{\struttext}{$\scriptstyle\mathstrut$}
\addtolength{\across}{.2cm}   
\def\row##1{\framebox[\across][l]{\boxify ##1, \closerow, }}
\def\vrow##1##2##3{\makebox[##1\sqlength][l]{}
                     \framebox[##2\sqlength][l]{\boxify ##3, \closerow, }}
\def\irow##1##2##3{\makebox[##1\sqlength][l]{}
                     \makebox[##2\sqlength][l]{\boxify ##3, \closerow, }}
\setlength{\vstrutht}{8.5pt} \addtolength{\vstrutht}{1.5\fboxsep} 
\setlength{\vstrutdp}{3.5pt} \addtolength{\vstrutdp}{1.5\fboxsep}
\setlength{\vstrutwd}{1pt} 
\else\ifx\testarg\compargw
\renewcommand{\struttext}{$\scriptstyle\mathstrut$}
\addtolength{\across}{.2cm}   
\def\row##1{\makebox[\across][l]{\boxify ##1, \closerow, }}
\def\vrow##1##2##3{\makebox[##1\sqlength][l]{}
                     \makebox[##2\sqlength][l]{\boxify ##3, \closerow, }}
\def\irow##1##2##3{\makebox[##1\sqlength][l]{}
                     \makebox[##2\sqlength][l]{\boxify ##3, \closerow, }}
\setlength{\vstrutht}{8.5pt} \addtolength{\vstrutht}{0\fboxsep} 
\setlength{\vstrutdp}{3.5pt} \addtolength{\vstrutdp}{0\fboxsep}
\setlength{\vstrutwd}{1pt} 
\else\ifx\testarg\compargt
\renewcommand{\struttext}{$\scriptscriptstyle\mathstrut$}
\addtolength{\across}{.2cm}   
\def\row##1{\framebox[\across][l]{\boxify ##1, \closerow, }}
\def\vrow##1##2##3{\makebox[##1\sqlength][l]{}
                     \framebox[##2\sqlength][l]{\boxify ##3, \closerow, }}
\def\irow##1##2##3{\makebox[##1\sqlength][l]{}
                     \makebox[##2\sqlength][l]{\boxify ##3, \closerow, }}
\setlength{\vstrutht}{8.5pt} \addtolength{\vstrutht}{1.5\fboxsep} 
\setlength{\vstrutdp}{3.5pt} \addtolength{\vstrutdp}{1.5\fboxsep}
\setlength{\vstrutwd}{1pt} 
\else
\def\row##1{\makebox[\across][l]{\boxify ##1, \closerow, }}
\def\vrow##1##2##3{\makebox[##1\sqlength][l]{}
                     \makebox[##2\sqlength][l]{\boxify ##3, \closerow, }}
\def\irow##1##2##3{\makebox[##1\sqlength][l]{}
                     \makebox[##2\sqlength][l]{\boxify ##3, \closerow, }}
\setlength{\vstrutht}{8.5pt} \addtolength{\vstrutht}{.5\fboxsep} 
\setlength{\vstrutdp}{3.5pt} \addtolength{\vstrutdp}{.5\fboxsep}
\setlength{\vstrutwd}{.5pt} 
\def\vbar{\llap{\makebox[4pt][l]{\smash{
\hbox{\vrule height14.0pt depth8.0pt width0.5pt}
}}}}
\fi\fi\fi\fi\fi
\begin{array}{l} }
\def\endtableau{\end{array}\endgroup}

\newcommand{\vstrut}{\hbox{\vrule height\vstrutht depth\vstrutdp width\vstrutwd}}
\newcommand{\vbar}{\llap{\makebox[2pt][l]{\smash{\vstrut}}}}
\newcommand{\row}[1]{\makebox[1in][l]{#1}}

\newdimen\passlength
\newdimen\passlengthtwo
\newbox\sizebox
\newlength{\sqlengthone}
\newlength{\sqlengthtwo}
\def\bitableau#1#2#3{\begingroup
\renewcommand{\arraystretch}{0}
\renewcommand{\vbar}{\smash{\vstrut}}
\settowidth{\sqlength}{$#1$}
\settowidth{\sqlengthone}{$#1$}
\settowidth{\sqlengthtwo}{$#2$}
\newcommand{\struttext}{\strut} 
\def\abox##1{\makebox[\sqlength][c]{\vphantom{\struttext} $##1\,$}}
\def\closerow{\relax}
\def\boxify##1,##2 {\def\testarg{##1}\def\comparg{\closerow}
        \ifx\comparg\testarg\let\next=\relax\else\abox{##1}\phantom{##2}
        \let\next=\boxify\fi\next}
\def\sboxify##1,##2 {\def\testarg{##1}\def\comparg{\closerow}
\ifx\comparg\testarg\let\next=\relax\else\mbox{
\vphantom{\struttext} $\,##1$}\,\phantom{##2}
\let\next=\sboxify\fi\next}
\def\bboxify##1,##2 {
\def\testarg{##1}\def\comparg{\closerow}
\ifx\comparg\testarg\let\next=\relax\else{\setbox\sizebox=\hbox{$\,##1$}
\ifdim\wd\sizebox>\passlength\global\passlength=\wd\sizebox
\fi}
\let\next=\bboxify\fi\next}
\def\bboxifytwo##1,##2 {
\def\testarg{##1}\def\comparg{\closerow}
\ifx\comparg\testarg\let\next=\relax\else{\setbox\sizebox=\hbox{$\,##1$}
\ifdim\wd\sizebox>\passlengthtwo\global\passlengthtwo=\wd\sizebox
\fi}
\let\next=\bboxifytwo\fi\next}
\def\comparg{f}  
\def\compargc{c} 
\def\compargs{s} 
\def\compargt{t} 
\def\compargw{w} 
\def\compargb{b} 
\def\testarg{#3}
\ifx\testarg\comparg
\addtolength{\across}{.2cm}   
\def\srow##1{\fbox{\sboxify ##1, \closerow, }}
\def\row##1{\setlength{\sqlength}{\sqlengthone}
                     \fbox{\boxify ##1, \closerow, }}
\def\rowt##1{\setlength{\sqlength}{\sqlengthtwo}
                     \fbox{\boxify ##1, \closerow, }}
\def\irow##1{\setlength{\sqlength}{\sqlengthone}
                     \hspace*{\fboxrule}\hspace*{\fboxsep}
                     \mbox{\boxify ##1, \closerow, }
                     \hspace*{\fboxrule}\hspace*{\fboxsep}}
\def\irowt##1{\setlength{\sqlength}{\sqlengthtwo}
                     \mbox{\boxify ##1, \closerow, }}
\setlength{\vstrutht}{8.5pt} \addtolength{\vstrutht}{1.5\fboxsep} 
\setlength{\vstrutdp}{3.5pt} \addtolength{\vstrutdp}{1.5\fboxsep}
\setlength{\vstrutwd}{1pt} 
\begin{array}{ll}
\else\ifx\testarg\compargc
\addtolength{\across}{.2cm}   
\def\|{\hfill\vrule}
\def\srow##1{\fbox{\sboxify ##1, \closerow, }}
\def\row##1{\setlength{\sqlength}{\sqlengthone}
                     \frame{\boxify ##1, \closerow, }}
\def\rowt##1{\setlength{\sqlength}{\sqlengthtwo}
                     \frame{\boxify ##1, \closerow, }}
\def\irow##1{\setlength{\sqlength}{\sqlengthone}
                     \hspace*{\fboxrule}\hspace*{\fboxsep}
                     \mbox{\boxify ##1, \closerow, }
                     \hspace*{\fboxrule}\hspace*{\fboxsep}}
\def\irowt##1{\setlength{\sqlength}{\sqlengthtwo}
                     \mbox{\boxify ##1, \closerow, }}
\setlength{\vstrutht}{8.5pt} \addtolength{\vstrutht}{1.5\fboxsep} 
\setlength{\vstrutdp}{3.5pt} \addtolength{\vstrutdp}{1.5\fboxsep}
\setlength{\vstrutwd}{1pt} 
\begin{array}{ll}
\else\ifx\testarg\compargs
\renewcommand{\struttext}{$\scriptstyle\mathstrut$}
\addtolength{\across}{.2cm}   
\def\srow##1{\fbox{\sboxify ##1, \closerow, }}
\def\row##1{\setlength{\sqlength}{\sqlengthone}
                     \fbox{\boxify ##1, \closerow, }}
\def\rowt##1{\setlength{\sqlength}{\sqlengthtwo}
                     \fbox{\boxify ##1, \closerow, }}
\def\irow##1{\setlength{\sqlength}{\sqlengthone}
                     \hspace*{\fboxrule}\hspace*{\fboxsep}
                     \mbox{\boxify ##1, \closerow, }
                     \hspace*{\fboxrule}\hspace*{\fboxsep}}
\def\irowt##1{\setlength{\sqlength}{\sqlengthtwo}
                     \hspace*{\fboxrule}\hspace*{\fboxsep}
                     \mbox{\boxify ##1, \closerow, }
                     \hspace*{\fboxrule}\hspace*{\fboxsep}}
\setlength{\vstrutht}{8.5pt} \addtolength{\vstrutht}{1.5\fboxsep} 
\setlength{\vstrutdp}{3.5pt} \addtolength{\vstrutdp}{1.5\fboxsep}
\setlength{\vstrutwd}{1pt} 
\begin{array}{ll}
\else\ifx\testarg\compargb
\def\srow##1{}
\def\row##1{{\bboxify ##1, \closerow, }}
\def\irow##1{}
\def\rowt##1{{\bboxifytwo ##1, \closerow, }}
\def\irowt##1{}
\begin{array}{ll}
\else\ifx\testarg\compargw
\renewcommand{\struttext}{$\scriptstyle\mathstrut$}
\addtolength{\across}{.2cm}   
\def\srow##1{\mbox{\sboxify ##1, \closerow, }}
\def\row##1{\setlength{\sqlength}{\sqlengthone}
                     \mbox{\boxify ##1, \closerow, }}
\def\rowt##1{\setlength{\sqlength}{\sqlengthtwo}
                     \mbox{\boxify ##1, \closerow, }}
\def\irow##1{\setlength{\sqlength}{\sqlengthone}
                     \mbox{\boxify ##1, \closerow, }}
\def\irowt##1{\setlength{\sqlength}{\sqlengthtwo}
                     \mbox{\boxify ##1, \closerow, }}
\setlength{\vstrutht}{8.5pt} \addtolength{\vstrutht}{0\fboxsep} 
\setlength{\vstrutdp}{3.5pt} \addtolength{\vstrutdp}{0\fboxsep}
\setlength{\vstrutwd}{1pt} 
\begin{array}{ll}
\else\ifx\testarg\compargt
\renewcommand{\struttext}{$\scriptscriptstyle\mathstrut$}
\addtolength{\across}{.2cm}   
\def\srow##1{\fbox{\sboxify ##1, \closerow, }}
\def\row##1{\setlength{\sqlength}{\sqlengthone}
                     \fbox{\boxify ##1, \closerow, }}
\def\rowt##1{\setlength{\sqlength}{\sqlengthtwo}
                     \fbox{\boxify ##1, \closerow, }}
\def\irow##1{\setlength{\sqlength}{\sqlengthone}
                     \hspace*{\fboxrule}\hspace*{\fboxsep}
                     \mbox{\boxify ##1, \closerow, }
                     \hspace*{\fboxrule}\hspace*{\fboxsep}}
\def\irowt##1{\setlength{\sqlength}{\sqlengthtwo}
                     \hspace*{\fboxrule}\hspace*{\fboxsep}
                     \mbox{\boxify ##1, \closerow, }
                     \hspace*{\fboxrule}\hspace*{\fboxsep}}
\setlength{\vstrutht}{8.5pt} \addtolength{\vstrutht}{1.5\fboxsep} 
\setlength{\vstrutdp}{3.5pt} \addtolength{\vstrutdp}{1.5\fboxsep}
\setlength{\vstrutwd}{1pt} 
\begin{array}{ll}
\else 
\def\srow##1{\mbox{\sboxify ##1, \closerow, }}
\def\row##1{\setlength{\sqlength}{\sqlengthone}
                     \mbox{\boxify ##1, \closerow, }}
\def\rowt##1{\setlength{\sqlength}{\sqlengthtwo}
                     \mbox{\boxify ##1, \closerow, }}
\def\irow##1{\setlength{\sqlength}{\sqlengthone}
                     \mbox{\boxify ##1, \closerow, }}
\def\irowt##1{\setlength{\sqlength}{\sqlengthtwo}
                     \mbox{\boxify ##1, \closerow, }}
\setlength{\vstrutht}{8.5pt} \addtolength{\vstrutht}{.5\fboxsep} 
\setlength{\vstrutdp}{3.5pt} \addtolength{\vstrutdp}{.5\fboxsep}
\setlength{\vstrutwd}{.5pt} 
\def\vbar{\llap{\makebox[4pt][l]{\smash{
\hbox{\vrule height14.0pt depth8.0pt width0.5pt}
}}}}
\begin{array}{ll}
\fi\fi\fi\fi\fi\fi
 }

\def\endbitableau{\end{array}\endgroup}

\newcommand{\bitblx}[2]{
\passlength=0pt
\passlengthtwo=0pt
\makebox[0pt][l]{$\bitableau{x}{x}b #2\endbitableau$}
\left(\bitableau{\hskip\passlength}{\hskip\passlengthtwo}#1 #2\endtableau
\right)}
\newcommand{\brbitblx}[2]{
\passlength=0pt
\passlengthtwo=0pt
\makebox[0pt][l]{$\bitableau{x}{x}b #2\endbitableau$}
\left[\bitableau{\hskip\passlength}{\hskip\passlengthtwo}#1 #2\endtableau
\right]}

\newcommand{\tightbrbitblx}[2]{
\passlength=0pt
\passlengthtwo=0pt
\makebox[0pt][l]{$\bitableau{x}{x}b #2\endbitableau$}
\left[\kern-8pt
\bitableau{\hskip\passlength}{\hskip\passlengthtwo}#1 #2\endtableau
\kern-6pt\right]}

\newcommand{\tblx}[2]{
\passlength=0pt
\passlengthtwo=0pt
\makebox[0pt][l]{$\bitableau{x}{x}b #2\endbitableau$}
\left[\bitableau{\hskip\passlength}{\hskip\passlengthtwo}#1 #2\endtableau
\right]}

\newcommand{\tblxshape}[2]{
\passlength=0pt
\passlengthtwo=0pt
\makebox[0pt][l]{$\bitableau{x}{x}b #2\endbitableau$}
\bitableau{\hskip\passlength}{\hskip\passlengthtwo}#1 #2\endtableau
}

\newcommand{\rsktblx}[2]{
\passlength=0pt
\passlengthtwo=0pt
\makebox[0pt][l]{$\bitableau{x}{x}b #2\endbitableau$}
\bitableau{\hskip\passlength}{\hskip\passlengthtwo}#1 #2\endtableau
}

\newcommand{\tshape}[1]{{  \hspace{-7pt}
  \newcommand{\trow}[1]{\row{##1}\phantom{{\smush{\row{\dims}}}}\\}
 {\setbox\strutbox=\hbox{\vrule height1pt depth1pt width0pt}
  \newcommand{\dims}{\vrule height3.75pt depth0pt width 3.80pt}
 {\newcommand{\x}{
\smash{\smush{\frame{\phantom{\vrule height4.40pt depth0pt width 5.20pt}}}}}
\tblxshape n{
#1
} } } }}

  \author{Brian D. Taylor\\
          Wayne State University\\
          Detroit, MI 48202, USA\\
          {\tt bdt@math.wayne.edu}}
  \title{A straightening algorithm 
         for row-convex tableaux.\thanks{\tt
         http://www.math.wayne.edu/\~\relax bdt/straightening.ps}}
  \date{July 26, 1999}
  
\begin{document}
\else
  \begingroup
  \let\done=\endgroup
\fi

\newcommand{\rsamesurround}[1]
{\begin{picture}(0,0)
\thicklines
\put(-25,3.75){\makebox(50,14.25)[t]{same}}
\put(0,0){\makebox[10pt][l]{#1}}
\put(-25,5.25){\oval(100,30)[r]}
\end{picture}
\hbox{\phantom{#1}}
}

\newcommand{\nonote}[1]{}

\newcommand{\SC}[1]{\hbox{\sc #1}}

\def\shflsign{{\rm sign}}

\newcommand{\wordtoinit}{\Psi}
\newcommand{\inittoword}{\Psi^{-1}}



\maketitle

\begin{abstract}
We produce a new basis for the
Schur and Weyl modules associated to a row-convex shape $D$.
The basis is indexed by new class of ``straight'' tableaux which we introduce
by weakening the usual requirements for standard tableaux. 
Spanning is proved via a new straightening
algorithm for expanding elements of the representation into this basis.
For skew shapes, this algorithm specializes to the classical straightening law.
The new straight basis is used to produce bases for flagged Schur and Weyl modules,
to provide Groebner and {\sc sagbi} bases for the homogeneous coordinate
rings of some configuration varieties and to produce a flagged branching rule for row-convex representations. 
Systematic use of supersymmetric letterplace techniques  
enables the representation theoretic results to be applied to
representations of the general linear Lie superalgebra as well as to the
general linear group.
\end{abstract}

\section{Introduction}
\nonote{rewrite this section and add some history. Include refs to initterms, thesis, applicn in Sot-Stu?, diff w. Br-Co,
applic-of-diag-term paper}
Akin, Buchsbaum, and Weyman in~\cite{ABW82} give a construction that associates 
a $GL_n$-representation to
any generalized { shape} like
\begin{equation}{  \hspace{-7pt}
 {\setbox\strutbox=\hbox{\vrule height1pt depth1pt width0pt}
  \newcommand{\dims}{\vrule height3.75pt depth0pt width 3.80pt}
 {\newcommand{\x}{
\smash{\smush{\frame{\phantom{\vrule height4.40pt depth0pt width 5.20pt}}}}}
\tblxshape n{
\row{  ,   , \x, \x }\phantom{{\smush{\row{\dims}}}}\\
\row{\x, \x, \x, \x }\phantom{{\smush{\row{\dims}}}}\\
\row{  , \x, \x     }\phantom{{\smush{\row{\dims}}}}\\
\row{\x, \x         }\phantom{{\smush{\row{\dims}}}}\\
} } } }
\quad \hbox{or}\quad
{  \hspace{-7pt}
 {\setbox\strutbox=\hbox{\vrule height1pt depth1pt width0pt}
  \newcommand{\dims}{\vrule height5pt depth0pt width 4.30pt}
 {\newcommand{\x}{
\smash{\smush{\frame{\phantom{\vrule height4.75pt depth0pt width 4.75pt}}}}
}
\tblxshape n{
\row{  , \x, \x     }\phantom{{\smush{\row{\dims}}}}\\
\row{\x,   , \x     }\phantom{{\smush{\row{\dims}}}}\\
\row{\x, \x         }\phantom{{\smush{\row{\dims}}}}\\
} } } }\quad \hbox{or}\quad
{  \hspace{-7pt}
 {\setbox\strutbox=\hbox{\vrule height1pt depth1pt width0pt}
  \newcommand{\dims}{\vrule height3.75pt depth0pt width 3.80pt}
 {\newcommand{\x}{
\smash{\smush{\frame{\phantom{\vrule height4.40pt depth0pt width 5.20pt}}}}}
\tblxshape n{
\row{  ,   , \x     }\phantom{{\smush{\row{\dims}}}}\\
\row{\x,   ,        }\phantom{{\smush{\row{\dims}}}}\\
\row{  , \x,        }\phantom{{\smush{\row{\dims}}}}\\
} } } }
\quad \hbox{or}\quad
{  \hspace{-7pt}
 {\setbox\strutbox=\hbox{\vrule height1pt depth1pt width0pt}
  \newcommand{\dims}{\vrule height3.75pt depth0pt width 3.80pt}
 {\newcommand{\x}{
\smash{\smush{\frame{\phantom{\vrule height4.40pt depth0pt width 5.20pt}}}}}
\tblxshape n{
\row{\x, \x, \x, \x }\phantom{{\smush{\row{\dims}}}}\\
\row{\x, \x,        }\phantom{{\smush{\row{\dims}}}}\\
\row{\x, \x         }\phantom{{\smush{\row{\dims}}}}\\
\row{\x,            }\phantom{{\smush{\row{\dims}}}}\\
} } } }
.\end{equation}
Significant progress has been made by Reiner and Shimozono and
by Lakshmibai and Magyar in describing bases for these $GL_n$-representations
and for the associated flagged representations of the
Borel subgroup of upper triangular matrices in $GL_n$. As one expects,
these bases are indexed by some subset of the generalized tableaux found
by filling each cell in the generalized shape with a number from 
$1$~to~$n$.

The present paper shows how to construct a well-behaved
{\em straight} basis for the representations
associated to any {\em row-convex} shape,
such as
{ $\hspace{-7pt}
 {\setbox\strutbox=\hbox{\vrule height1pt depth1pt width0pt}
  \newcommand{\dims}{\vrule height3.75pt depth0pt width 3.80pt}
 {\newcommand{\x}{
\smash{\smush{\frame{\phantom{\vrule height4.40pt depth0pt width 5.20pt}}}}}
\tblxshape n{
\row{  ,   , \x, \x }\phantom{{\smush{\row{\dims}}}}\\
\row{\x, \x, \x     }\phantom{{\smush{\row{\dims}}}}\\
\row{  , \x         }\phantom{{\smush{\row{\dims}}}}\\
\row{\x,            }\phantom{{\smush{\row{\dims}}}}\\
} } }$}, with no gaps in any row.
In particular, we give a 
local condition for testing whether a tableau is straight, we
give a straightening law that modifies only two rows at a time,
and the basis we present reduces immediately to a flagged basis.
The straight basis is distinct from the bases produced by Reiner-Shimozono
and Lakshmibai-Magyar, but see~\cite{initterms} for some combinatorial and
algebraic relationships between these bases. The straight bases provide a canonical
choice of basis for certain row-convex and column-convex ``almost-skew shapes.''
These shapes were shown by Woodcock in~\cite{W94} to possess a class of easily
flagged bases, but no method was presented for distinguishing a basis
in this class or for straightening elements of the representation into
a linear combination of basis elements.

Results on flagged tableaux are deduced in Section~\ref{sect:flagged}
from the main theorem on straight bases.
As shown in Section~\ref{sect:groebner}, the straight basis and straightening algorithm
may be applied to produce quadratic Groebner bases and  {\sc sagbi} bases for the homogeneous
coordinate rings of certain configuration varieties. Further applications
to commutative algebra may be found in~\cite{thesis} and~\cite{straightApplic}.
Applications to the representation theory of $GL_n$, $B_n$, $S_n$, and the general linear
Lie superalgebras are derived in Section~\ref{sect-br-rule}
where a branching rule is produced for decomposing a row-convex $GL_n$-representation
in terms of $GL_{n-1}$-representations.

This paper studies the Schur and Weyl modules as special cases of the
super Schur modules  which we construct as submodules of the 
letterplace superalgebra. All results in this paper are characteristic--free
and the requisite background on superalgebras is detailed in Section~\ref{SectSuper}.
Much of the presentation in Section~\ref{SectSuper} is new and, we hope, accessible
to the non-specialist. The construction proper is given in
Section~\ref{SectConstr}. Straight tableaux are introduced and
independence is proved in Section~\ref{SectBasis}. Section~\ref{sect:straighten},
the heart of the paper, presents the straightening algorithm.

\section{Polynomial superalgebras}
\label{SectSuper}

This section introduces the definitions required to make the main results
of this paper characteristic free and applicable to Weyl modules. The reader
concerned only with Schur modules 
in characteristic~0 may safely take $\C L$ and $\C P$ to be
the positive integers, $\N$, (or finite subsets of $\N$.) The set $\C L$ may be thought
of as the indexing the rows of a generic matrix  $(x_{i,j})$ and $\C P$ indexes the columns.
We may then take take $\Superlp LP$ to be the polynomial ring whose variables are matrix entries $x_{i,j}$.
The {\em letterplace} $(i|j)$ is taken to be shorthand for $x_{i,j}$ and the
expression $[i_1,\ldots,i_k|j_1,\ldots,j_k]$ is taken to be the determinant of the
$k\times k$ minor  $(x_{i_r,j_s})$ of the matrix $(x_{i,j})$.

The constructions used in this paper take place inside polynomial superalgebras over the
integers, $\Z$, that is inside tensor products of symmetric, exterior, and divided powers algebras.
We construct the polynomial superalgebras over $\Z$ as $\Z$-subalgebras of a symmetric algebra over
the rationals, $\Q$, tensored with an exterior algebra over $\Q$.
Write the symmetric and exterior $\integers$-algebras associated to a set $\C L$ 
as $\Sym(\C L)$\index{\Sym(\C L)} and $\Wedge(\C L)$\index{\Wedge(\C L)}.
These are $\Z$-subalgebras of the symmetric and exterior
$\rationals$-algebras 
 $\Sym_\Q(\C L)$\index{\Sym_\Q(\C L)}
and $\Wedge_\Q(\C L)$\index{\Wedge_\Q(\C L)} associated to $\C L$.
The divided powers algebra, $\Div(\C L)$ of a 
set $\C L$ is the $\Z$-subalgebra of $\Sym_\Q(\C L)$ generated by all $x^i\over i!$
for all $x\in\C L$.

We define a \defi{signed set} to be a set $\C L$ together with a
a function $|\,|:\C L\rightarrow\Z_2$. 
\label{ineqdef}We say that elements in the preimage of $0$
are \defi{positively signed}; we call this preimage $\C L^+$.
Elements in the preimage, $\C L^-$, of $1$ are said to be \defi{negatively signed}.
A signed set $\C L$ endowed with a total order, $<$, is called a (signed) {\em alphabet}.
For notational convenience, we define two new inequalities, $\lp$ and $\lm$ on $\C L$. We
say that $a\lp b$ (respectively $a\lm b$) when $a<b$ or when $a=b$ and $|a|=|b|=0$ (respectively
$a=b$ and $|a|=|b|=1$.) 

A {\em superalgebra} is simply an algebra with a $\Z_2$-grading. 
We construct a $\Q$-superalgebra with the elements of a signed set as generators and such that the grading
on these generators is $|\,|$.
For any signed set $\C L$,
define $\SuperQ{\C L}$ to be  $\Sym_\Q(\C L^+)\otimes\Wedge_\Q(\C L^-)$.
Likewise, we define $\Super{\C L}$ to be $\Div(\C L^+)\otimes\Wedge(\C L^-)$; as above, we may
consider this to be a $\Z$-subalgebra of $\SuperQ{\C L}$.
Given another signed set $\C P$, we will define the ``letter-place'' algebra, $\Superlp LP$, to be a
$\integers$-subalgebra of $\Superembed LP$ where $|x_{a,d}| = |a|+|d|$. 
In particular, $\Superlp LP$ is the subalgebra
generated by all $x_{a,d}$ and by all $x_{a,d}^i\over i!$ with $a,d\in\C L^+$,
$i\in\naturals$.
This algebra is naturally isomorphic to 
\[   \Wedge(\C L^-\times \C P^+ \dUnion \C L^+\times \C P^-)\tens 
     \Sym(\C L^-\times \C P^- )\tens 
     \Div(\C L^+\times \C P^+ ) 
\,.\]
We extend $|\,|$ to a $\Z_2$ grading of $\Superlp LP$.
Following~\cite{GRS}, we write the elements $x_{a,d}$ of $\Superlp LP$
as the signed variables $(a|d)$\index{(a|d)}, and we will
define the \defi{biproduct}, $(w_1,\ldots,w_k|v_1,\ldots,v_k)$, of a pair of 
sequences $\mbf w$ and $\mbf v$ in $\C L$ and $\C P$ respectively.

\begin{definition}\label{signrule}
Given  sequences $\mbf w=w_1,\ldots, w_k \in \C L$
and              $\mbf v=v_1,\ldots, v_k \in \C P$,
define
\index{(\mbf w | \mbf v)}
$
(\mbf w | \mbf v) =(w_1,\ldots,w_k|v_1,\ldots,v_k)=
\sum_{\sigma\in S_k} (-1)^{n_\sigma} 
   (w_{\sigma(1)}|v_1)  \cdots (w_{\sigma(k)}|v_k) 
$
\[\displaylines{where\quad
n_\sigma = {\# \{ (i,j) : i<j,\ \sigma^{-1}(i) > \sigma^{-1}(j),\ 
       \hbox{$w_i,w_j$ are negative} }\}
\cr\hfill
+{\# \{ (i,j) : i>j\ \hbox{and $w_{\sigma(i)}, v_j$ are negative} \} }.
\qquad}\]
\end{definition}
The following definition/proposition  indicates that the
biproduct can be thought of as a bilinear map on $\Super{\C L}\times\Super{\C P}$.

\begin{definition}\label{bilinear}
Given  sequences $w=w_1,\ldots, w_k \in \C L$
and              $v=v_1,\ldots, v_k \in \C P$,
define
$(w_1\cdot w_2 \cdots w_k|v_1\cdot v_2 \cdots  v_k) =
(w_1, w_2, \ldots, w_k|v_1, v_2, \ldots,  v_k)$. Extend this by bilinearity
to a map $(\ |\ ):\Super{L}\times\Super{P}\rightarrow\Superlp LP$.
\end{definition}
It is straightforward to to check that this map is well-defined.

In order to better handle divided powers as elements of a $\Z$-subalgebra 
contained in a $\Q$-algebra, we make the following definition.
If $w$ is a sequence in $\C A$ for some signed set $\C A$,
then we define $\c(w)!$
\index{\c(w)!}
to be $\prod_{i\in\C A^+} (\#\hbox{times $i$ appears in $w$})!$.

Call a monomial $M=\prod_i(l_i|p_i)\in \SuperQlp LP$ {\em sorted} when 
$p_1\le p_2\le \cdots$ and $p_i=p_{i+1}\in\C P^+$ implies $l_i\lp l_{i+1}$
and dually $p_i=p_{i+1}\in\C P^-$ implies $l_i\gm l_{i+1}$.
$\Superlp LP$ is a free $\Z$-module with basis consisting
of the {\em divided powers monomials}
$\{ {1\over \c(M)!} M \}$ for all sorted monomials $M\in\SuperQlp LP$.
Here $\c(M)!=\c({(l_1|p_1),(l_2|p_2),\ldots})!$.
We will consider two monomials (respectively divided powers monomials)
the same when they differ by a nonzero scalar multiple
(respectively a multiple of $\pm1$.)

Define a function
         $\Tab(\mbf w| \mbf v)$
when $\mbf w$ (respectively $\mbf v$) is a $k$-tuple of letters 
in $\C L$ (respectively $\C P$) by
\index{ \Tab(w_1,\ldots, w_k | v_1,\ldots, w_k) }
\[ \Tab(w_1,\ldots, w_k | v_1,\ldots, w_k) =
{1\over \c(\mbf w)! \c(\mbf v)!} 
    (-1)^{\#\,\{(i,j) : i>j,\ 
         w_i\in\C L^-, v_j\in\C P^+ \} }
    (w|v).
\]
Observe that the divided powers monomials occur with coefficient $\pm1$ in
the expansion of $\Tab(\mbf w|\mbf v)$ and that if
$w_1\lp w_2\lp \cdots\lp w_k$ and $v_1\lp v_2\lp \cdots \lp v_k$,
then the basis element  ${1\over \c({\prod_i (w_i|v_i)} )!} \prod_i (w_i|v_i)$
appears with coefficient~$1$.




\section{Schur modules, Weyl modules, generalizations}
\label{superSchur}
\label{SectConstr}

In this section, we define our primary object of study,
the super-Schur module as a $\integers$-submodule of a letterplace algebra.
The unsigned cases produce the Schur and Weyl modules of Akin-Buchsbaum-Weyman~\cite{ABW82}
when $\C L=\C L^-$ and $\C L=\C L^+$ respectively.

We define two tableaux associated to a shape. The first
is useful for referring to cells in the shape and the second
plays a fundamental role in our construction of the super-Schur
modules. 

\begin{definition}
Let $D$ be a shape. 

Define $F(D)$ to be the tableau 
of shape $D$ whose cells are labeled $1,2,3,\ldots$
starting with the northmost cell in the leftmost column
and continuing down the column, then down the second
leftmost column, etc. In this paper, the signs of the
letters in $F(D)$ are irrelevant.

A tableau of shape $D$ is termed \defi{Deruyts} if it is
 obtained by filling each
cell in the diagram with the cell's column index viewed as a negative
variable. We denote such a tableau by $\Deruyts(D)$.
\index{\Deruyts(D)}

\end{definition}

{\em Shapes appearing in this paper are assumed, unless otherwise noted, 
to have first coordinate~$1$ in their top rows
and  second coordinate~$1$ in their leftmost columns. }

\begin{example}
\[ F\left( \tshape{
\trow{  ,   , \x, \x }
\trow{\x,   , \x, \x }
\trow{  ,   , \x     }
\trow{\x,            }
}\right) =
\rsktblx w{
\row{    ,   ,  3  ,  6   }\\
\row{ 1  ,   ,  4  ,  7   }\\
\row{    ,   ,  5         }\\
\row{ 2  ,                }\\
}
\qquad \qquad
\Deruyts\left( \tshape{
\trow{  ,   , \x, \x }
\trow{\x,   , \x, \x }
\trow{  ,   , \x     }
\trow{\x,            }
}\right) =
\rsktblx w{
\row{    ,   ,  3^-,  4^- }\\
\row{ 1^-,   ,  3^-,  4^- }\\
\row{    ,   ,  3^-       }\\
\row{ 1^-,                }\\
}
\]
\end{example}

\begin{definition}
Suppose $S$ and $T$ are tableaux of the same shape.
Let the word ${\SC s}_i$ be the $i$th row of $S$ and let ${\SC t}_i$ be the $i$th row of $T$.
Define ${[S|T]=\prod_i\Tab({\SC s}_i|{\SC t}_i)}$.\index{[S|T]}
Hence $[S|T]=\prod_i [{\SC s}_i|{\SC t}_i]$.

Now suppose that $T$ is a tableau of shape $D$.
Suppose that $\C L$ contains the set of letters present in $T$ and that $\C P^-$
contains the indices for all columns present in $D$.
Define an element
\index{[T]}
$[T]\in \Superlp LP$ indexed by $T$ by,
\[ [T]  =  \Big[ T \Bigm| \Deruyts(D) \Big] .\] 
\end{definition}
\begin{example}
Let ${\cal L}={\cal L}^- =\{a,b,c,d,e,f,g\}$ and let 
${\cal P}={\cal P}^-=\{1,2,3,4\}$. 
Let $T=\rsktblx w
              {\row{ , , a, d}\\
               \row{b, , c, e}\\
               \row{ , , f,  }\\
               \row{g, ,  ,  }}
$.
Then 
\[ [T] =   \brbitblx w{\row{ , , a, d} &\vbar\ \rowt{ , , 3, 4}\\
                       \row{b, , c, e} &\vbar\ \rowt{1, , 3, 4}\\
                       \row{ , , f,  } &\vbar\ \rowt{ , , 3,  }\\
                       \row{g, ,  ,  } &\vbar\ \rowt{1, ,  ,  }   }
          .\]
In other words, $[T]$ is a scalar multiple of  $(ad|23)(bce|123)(f|2)(g|1)$.
The scalar in this example being $1$, we have
\[ [T] = 
\det\left( \begin{array}{cc}
(d|3) & (d|4) \\
(a|3) & (a|4) \\
\end{array}\right)
\det\left( \begin{array}{ccc}
(e|1) & (e|3) & (e|4) \\
(c|1) & (c|3) & (c|4) \\
(b|1) & (b|3) & (b|4) \\
\end{array}\right)
(f|2)
(g|1)\,\,.\]
\end{example}

\begin{definition}
Suppose that $D$ is a shape. Define
the \defi{super-Schur module}
\index{\tabmod_{\C D}(\C L)}
\[ \tabmod_{D}(\C L) = {\rm span}_\integers 
           \Big\{ [T] : {\rm shape}(T) = {\cal D} \hbox{ and $T$ is filled
                         with letters from $\C L$} \Big\} .\] 
\end{definition}
In the case that $\cal L$ is negative (respectively
positive) then $\tabmod_{D}(\C L)$ is called the Schur (respectively Weyl) module
associated with the diagram $D$. These terms are justified by the following result.
A proof may be found in~\cite{thesis}.

\begin{proposition}\label{asComposite}
Let $R$ be a commutative ring. Let $F$ be a free $R$-module of rank $n$.
Let $\alpha$ be the $0/1$-matrix having $1$'s precisely where $D$ has cells.

If $\C L=\C L^-$ has cardinality $n$, then
$R\tensZ\tabmod_D(\C L^-) = L_\alpha(F)$, where $L_\alpha(F)$ is the Akin-Buchsbaum-Weyman
``Schur functor'' associated to the generalized shape matrix $\alpha$.

If $\C L=\C L^+$ has cardinality $n$, then
$R\tensZ\tabmod_D(\C L^+) = K_\alpha(F)$, where $K_\alpha(F)$ is the Akin-Buchsbaum-Weyman
``coSchur functor.''
\qed\end{proposition}

\begin{example}\label{exmplbasis}
The Weyl module of shape 
${  \hspace{-7pt}
 {\setbox\strutbox=\hbox{\vrule height1pt depth1pt width0pt}
  \newcommand{\dims}{\vrule height3.75pt depth0pt width 3.80pt}
 {\newcommand{\x}{
\smash{\smush{\frame{\phantom{\vrule height4.40pt depth0pt width 5.20pt}}}}}
\tblxshape n{
\row{\x, \x, \x,  }\phantom{{\smush{\row{\dims}}}}\\
\row{  , \x,      }\phantom{{\smush{\row{\dims}}}}\\
} } } }$ 
on positive letters $a,b$ is spanned by
{\small \[\displaylines{
\tblx w{
\row{a, a, a}\\
\row{ , b,  }\\
},\quad
\tblx w{
\row{a, a, b}\\
\row{ , b,  }\\
},\quad
\tblx w{
\row{b, b, b}\\
\row{ , a,  }\\
}
,\quad
\tblx w{
\row{a, b, b}\\
\row{ , a,  }\\
}
,\quad
\tblx w{
\row{a, b, b}\\
\row{ , b,  }\\
}
,\quad
\tblx w{
\row{a, a, b}\\
\row{ , a,  }\\
}
,\cr
\tblx w{
\row{a, a, a}\\
\row{ , a,  }\\
}
,\quad
\tblx w{
\row{b, b, b}\\
\row{ , b,  }\\
}
}\] }
in the exterior algebra generated by the anti-commuting variables
$(a|1)$, $(a|2)$, $(a|3)$, $(b|1)$, $(b|2)$, $(b|3)$.
The last two of the above skew-polynomials are identically $0$. In the
next section we single out the first three elements as a basis.
\end{example}

\section{Row-convex diagrams and straight tableaux}
\label{SectBasis}

The usual bases for skew Weyl modules consist of the
semistandard Young tableaux, namely all tableaux which
weakly increase in their rows and strictly increase in their columns.
Example~\ref{exmplbasis} showed that this is not
the case for more general shapes. 
Nevertheless, the basis of~\cite{ABW82} for shape $D$ skew Weyl
modules indexed by standard Young tableaux of shape $D$
has a number of properties we wish
to preserve. In particular:
\label{niceprops}\begin{enumerate}
\item \label{incrprop}
The rows of the tableaux in the indexing set weakly increase.
\item \label{diagprop}
Knowing the number of times a letter appears
in each column of a tableau in the indexing set determines
that tableau.
\item It is combinatorially ``obvious'' when a tableau
is in the indexing set.
\item The elements $[T]$ where $T$ is in the index set form
a basis for module. \label{proddets}
\item There is an easy to describe algorithm for
rewriting $[T]$ in terms of basis elements.
\end{enumerate}
Property~\ref{diagprop} underlies the {\sc sagbi}-basis
algorithms of~\cite{Stu93}; in~\cite{W94} Woodcock
shows that there must exist bases satisfying this
property when $D$ is ``almost-skew.''

Only slightly more complicated shapes,
$\tshape {
\trow{\x,   , \x }\\
\trow{\x, \x,   }\\}$
for instance, fail to simultaneously possess properties~\ref{incrprop}, \ref{diagprop}, and~\ref{proddets}.
To see this, examine the Specht module associated to this shape. Recall that this is the subspace of the associated 
Schur module spanned by all tableaux containing letters $1^-,2^-,3^-,4^-$ with no repeats; here the indexing shape is
transposed from the indexing shape used in~\cite{Sa}.
This Specht module is isomorphic to the one indexed by  
$\tshape {
\trow{  , \x, \x }\\
\trow{\x, \x,   }\\}$
hence has dimension~5. However, there are only 4 tableaux of shape 
$\tshape {
\trow{\x,   , \x }\\
\trow{\x, \x,   }\\}$
satisfying conditions~ \ref{incrprop} and \ref{diagprop}.

We define a class of ``straight'' tableaux 
satisfying the above properties.
The elements $[T]$ where $T$ is straight and of shape $D$ will form a basis
for the super-Schur module $\tabmod_D$ for any ``row-convex''
shape $D$.

\begin{definition}
A \defi{row-convex shape},  such as
{ $\hspace{-7pt}
 {\setbox\strutbox=\hbox{\vrule height1pt depth1pt width0pt}
  \newcommand{\dims}{\vrule height3.75pt depth0pt width 3.80pt}
 {\newcommand{\x}{
\smash{\smush{\frame{\phantom{\vrule height4.40pt depth0pt width 5.20pt}}}}}
\tblxshape n{
\row{  ,   , \x, \x }\phantom{{\smush{\row{\dims}}}}\\
\row{\x, \x, \x     }\phantom{{\smush{\row{\dims}}}}\\
\row{  , \x         }\phantom{{\smush{\row{\dims}}}}\\
\row{\x,            }\phantom{{\smush{\row{\dims}}}}\\
} } }$}, is  a shape with no gaps in any row.  
I.E., if cells $(r,i)$ and $(r,k)$ are in a shape $D$,
then $(r,j)$ is in $D$, for all $i<j<k$.
Since the constructions of section~\ref{superSchur}
 are not sensitive to the order of rows in a diagram, 
we assume that {\em the rows of a row
convex diagram are sorted so that higher rows end at least as far to the
right as lower rows}.
\end{definition}
We can denote any row-convex shape by $\lambda/\mbf m$ where $\lambda$ is a partition
and $\mbf m$ is a composition satisfying $m_i\le \lambda_i$ for all $i$; a cell is
in position $(i,j)$ of $\lambda/\mbf m$ iff $m_i\le j\le \lambda_i$.

Following~\cite{GRS}, and employing the notation for inequalities introduced
on page~\pageref{ineqdef}, a tableau $T$ with entries in a signed set
is {\em standard}\glossary{standard tableau} when it 
($\lp$)-increases across rows and $(\lm)$-increases down columns.

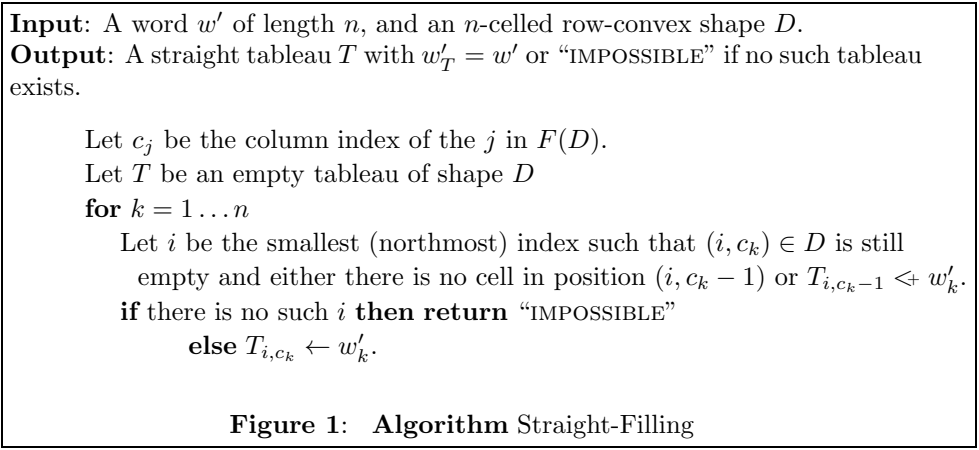
\begin{figure}[!tb]
\framebox{\vbox{
{\bf Input}: A word $w'$ of length $n$, and an $n$-celled row-convex shape $D$.
\\
{\bf Output}: A straight tableau $T$ with $w'_T=w'$ or 
      ``{\sc impossible}'' if no such tableau exists.
\begin{code}
Let $c_j$ be the column index of the $j$ in $F(D)$. \\
Let $T$ be an empty tableau of shape $D$\\
\For $k=1\ldots n$\\
\>Let $i$ be the smallest (northmost) index such that 
$(i,c_k)\in D$ is still\\
\>\ \ empty and either there is no cell
    in position $(i,c_k-1)$  or $T_{i,c_k-1}\lp w'_k$.\\
\>\If there is no such $i$ {\bf then} \Return ``{\sc impossible}''\\
\>\>{\bf else} $T_{i, c_k}\gets w'_k$.
\end{code}
\caption[Algorithm Straight-Filling]{
{\bf Algorithm} Straight-Filling \label{strt-fill:fillByLetter}}
}}
\end{figure}

I introduce the notion of a \defi{straight tableau} of row-convex shape
by slightly relaxing the usual conditions for standardness of a  tableau.
\begin{definition}\label{straight}
A row-convex tableau is called {\em straight} when 
\begin{enumerate}
\item \label{condone}
The contents of any row 
$\lp$-increase from left to right, and
\item \label{condtwo}
Given two cells in the same column, say $(i,k)$ and $(j,k)$ for
$i<j$, the entry in the top cell, $(i,k)$, may be $(\gp)$-larger than the entry in 
$(j,k)$ 
(i.e. the cells form an \defi{inversion})
only if cell $(i,k-1)$ exists and its content is $(\gm)$-larger than the
content of $(j,k)$.
\end{enumerate}
\end{definition}
This definition amounts to requiring that the columns are as close as possible to 
($\lm$)-increasing, subject to the condition that the rows remain ($\lp$)-increasing.
A more precise version of the preceding fact is implicit in the correctness
Algorithm Straight-Filling in Figure~\ref{strt-fill:fillByLetter}.
A tableau satisfying condition~\ref{condone} is 
called \defi{row-standard} and an inversion
violating condition~\ref{condtwo} is called a {\em flippable} inversion.

\begin{proposition}
A skew tableau, $T$, is straight iff it is standard.
\end{proposition}
\begin{proof}
Since a standard tableau has no inversions, it
suffices to prove the only-if part. We prove the contrapositive.
We can assume that $T$ is row-standard.
Suppose  that
the cells $(i,k),(j,k)$ with $i<j$ are an inversion. Let $k_0$
be the least (leftmost) column such that $(i,k_0),(j,k_0)$ is
an inversion. 
If $(i,k_0-1)$ exists then by skewness so does $(j,k_0-1)$ and thus
by assumption $T_{i,k_0-1}\lm T_{j,k_0-1} \lp T_{j,k_0}$ hence
$T$ is not straight.
\end{proof}
\begin{corollary}
The straight tableaux of skew shape with only 
positively signed letters are the usual semistandard Young tableaux. \qed
\end{corollary}

\nonote{A word of intro before the defn?}
\begin{definition}
Given a tableau $T$ its \defi{column word}, $c_T$\index{u_T}
is the word formed by reading the entries of $T$ from bottom to top and
left to right.
Its \defi{modified column word} is the word $w_T$\index{w_T}
formed by writing the entries of the first column in weakly decreasing order
followed by the entries of the second column in decreasing order, etc.

We shall also occasionally require a \defi{reverse column word} $w'_T$\index{w'_T}
of $T$ formed by writing the entries of the first column of $T$ in increasing order
then those of the second column in increasing order etc.
\end{definition}

\begin{theorem}\label{uniqColCont}
If $T$ and $T'$ are straight tableaux of the same shape, then $T\ne T'$
implies $w_T\ne w_{T'}$.
More strongly, if there exists a straight tableau $T$ of shape $D$ with
$w_T=w$ then the algorithm {\tt Straight-Filling} in 
Figure~\ref{strt-fill:fillByLetter} produces it.
\end{theorem}

\nonote{do an example?}
\begin{proof}
A tableau, $T$, produced by this algorithm
must be straight. If in a fixed column, $k$, 
the letter $y$ is inserted into row $i$ by the algorithm
while $x\lp y$ was inserted into row $j>i$, then it
must be that $T_{i,k-1}\gm x$ else cell $T_{i,k}$ would
have been available to $x$ hence $x$ would have been 
placed there.

Now suppose that the algorithm produces a tableau $T$ with
reverse column word $w'$.
Let $\mbf c$ be as in the algorithm. 
Any tableau with reverse column word $w'$ can be
produced by a similar filling process. 
Define $\mbf i$ so that
reading through $w'$ and
inserting $w'_k$ into cell $(i_k,c_k)$ gives the desired
tableau. Let us assume that
if $w'_k$ appears in multiple cells in column $c_k$ that
the first $w'_k$ in $w'$ is used to fill the northmost appearance in the column, 
the second is used to fill the second northmost appearance, etc.

Let $\mbf i$ be the filling sequence corresponding to $T$,
this is the sequence produced by the Algorithm~Straight-Filling.
Let $\mbf i'$ be the filling sequence corresponding to some other
tableau $T'$.
Let $k_0$ be the smallest integer such that $i_{k_0}\ne i'_{k_0}$.
So in filling $T'$,
we have placed $w'_{k_0}$ into cell $(i'_{k_0},c_{k_0})$
when according to Algorithm~Straight-Filling,
it could have been put into $(i_{k_0},c_{k_0})$
where $i_{k_0}<i'_{k_0}$. By necessity, in filling $T'$,
 something $(\ge)$-larger than 
$w_{k_0}$ must be placed in $(i_{k_0},c_{k_0})$.
By our  assumptions about repeated letters in
the definition of $\mbf i$, this inequality is strict.
But these facts guarantee that the inversion $\{
(i_{k_0},c_{k_0}),(i'_{k_0}, c_{k_0})\}$ of $T'$ violates
condition~\ref{condtwo} in the definition of straight tableaux.

The above argument says that if we 
try to create a straight tableau $T$ with $w'_T=w'$
by reading across $w'$ 
and sequentially filling its letters into a tableau
then at each step the choice of where to insert the letters
is forced on us. If at any point during execution of the algorithm
there is no place to put a letter which preserves row-standardness,
then it is in fact not possible to find a straight tableau with
the designated column content and shape. This is precisely the
circumstance under which ``{\sc impossible}'' is returned.
  
We conclude that not only does {\tt Straight-Filling} produce
a straight tableau, but any other tableau, $T'$ having the
same modified (equivalently reverse) column word is not straight.
\end{proof}

\begin{corollary}\label{indep}
The matrix expressing the super-polynomials $[T]$ 
indexed by straight tableaux as $\Z$-linear
combinations of divided powers
monomials in the polynomial superalgebra is in echelon form with
$\pm 1$ at each pivot.
Hence the straight basis elements are linearly independent.
\end{corollary}

We defer the proof in order to develop the appropriate orders on basis
elements and monomials.
Monomials are ordered according to a generalization of the ``diagonal term order''
in~\cite{Stu93} which requires that the smallest monomial in $\det(A)$,
where $A$ is a minor of $(x_{i,j})$, be the product of the elements on the
diagonal. 
For compatibility with lexicographic order in Lemma~\ref{normalized-monoms}
this is backwards from the convention  in commutative algebra
which has $\prod_i(x_{i,i})$ be the {\em largest} monomial in $\det(A)$.

\begin{definition}
A \defi{diagonal term order} on $\Superlp LP$ is
\begin{enumerate}
\label{diagdefs}
\item  A total order, $\prec$, on monomials in $\Superlp LP$
       such that for monomials $m,m',n,n'$, the relations
       $m\prec m'$ and
       $n\prec n'$ imply that $mn\prec m'n'$ or $mn=0$ or
       $m'n'=0$.

\item The smallest monomial in a nonzero biproduct
      $(i_1\,\ldots,i_k | j_1,\ldots, j_k)$ with $i_1<\cdots< i_k$ 
      and $j_1<\cdots< j_k$ is $\prod_l (i_l|j_l)$.
\end{enumerate}
\end{definition}

The \defi{default diagonal term order}, $\diag$\index{diag}
that we utilize is characterized below.
\nonote{``reverse lexicographic'' term order in comm. alg. (double check)}
We order letterplaces $(i|j)$ by $(i|j) > (i'|j')$
when $j<j'$ or when $j=j'$ and $i>i'$. 
Let $M,N\in\Superlp LP$ be two
nonzero monomials. Suppose $(i|j)$ is the largest letterplace appearing
to a different power in $M$ and $N$. 
Write $N\diag M$ when $M$ is divisible by a higher power of $(i|j)$ than is $N$.
\begin{example}
Suppose $\C L = \{1^-, 2^-\}$ and $\C P =\{ a^+, x^- \}$,
then $(2|a) > (1|a) > (2|x) > (1|x)$. Further, we have
\[ (1|x)<(1|x)^6<(2|a)(1|x)^2(2|x)<(2|a)(2|x)^2(1|x)<(2|a)(1|a)(1|x)^2 .\]
\end{example}

The following lemma is immediate.
\begin{lemma}\label{normalized-monoms}
A \defi{normalized monomial} $\prod_{l=1}^k (i_l|j_l)\ne0$ in $\Superlp LP$ 
is a monomial written so that $(i_l|j_l)\ge (i_{l+1}|j_{l+1})$ in
the default diagonal term order.
For two normalized monomials,
$M=\prod_{l=1}^k (i_l|j_l)$ and $N=\prod_{l=1}^k (i'_l|j_l)$ 
differing only in their letters, $M<N$
in the default diagonal term order  iff
$i_1,\ldots,i_k$ is lexicographically less than $i'_1,\ldots, i'_k$.
\qed
\end{lemma}

\begin{definition}
Let \indx{\Psi} be the function taking a normalized monomial 
      $\prod_{l=1}^k (i_l|j_l)\in\Superlp LP$ 
to $i_1,\ldots,i_k$.
\end{definition}
 
\begin{definition}
Given $p\in\Superlp LP$ and an order $\prec$ on monomials, define
the \defi{initial monomial}
\indx{\inwt\prec{p}} of $p$ to be the smallest 
(divided powers) monomial appearing in $p$.

Sometimes the phrase ``initial term'' will be used when the coefficient
of the initial monomial is to be included.
\end{definition}
The following result says that in most cases the modified
column word of $T$ can be read directly from the smallest monomial
appearing in $[T]$.
\begin{proposition}
\label{onesixsix}
If $T$ is a  tableau whose rows ($\lp$)-increase
and whose columns contain no repeated positive letters, then
\[  w_T = \Psi(\inwt{\diag}{[T]}) .\] 
\end{proposition}
\begin{proof}
Suppose
$ [T]=\prod_i 
\Tab(w_{i,c_{i,1}},w_{i,c_{i,2}},\ldots,w_{i,c_{i,l_i}} |
               c_{i,1}, c_{i,2}, \ldots, c_{i,l_i} ) .$
The initial term (with coefficient) of the $i$th multiplicand is
$         \prod_j (w_{i,c_{i,j}} | c_{i,j} )   .$
and since positive letters never repeat in a column the product
of these initial terms is nonzero and hence equals $\inwt{\diag}{[T]}$.
\end{proof}

Note that the initial term
$     \prod_i  \prod_j (w_{i,c_{i,j}} | c_{i,j} )   $
appearing above is (up to sign) a basis element in the monomial $\Z$-basis for
$\Superlp LP$ and we have proved the following.

\begin{proposition}\label{1.5.6p}
If $T$ is straight of shape $D$, and if $c_l$ is the index
of the column of $F(D)$ containing $l$ then 
\[  \prod_l \Big( {(w'_T)}_l \Big| c_l \Big) =
\inwt{\diag}{[T]}\] 
The coefficient of the initial monomial of $[T]$ 
is $\pm 1$.\qed 
\end{proposition}

\begin{corollary}\label{readColCont}
Suppose  $T$ is a  straight tableau, then
$  w_T = \Psi (\inwt{\diag}{[T]}).$ \qed
\end{corollary}

We now complete the proof of the independence result.
\begin{proof}(of Corollary~\ref{indep}.) 
Since Theorem~\ref{uniqColCont} says that distinct straight tableaux
have distinct modified column words, we conclude from
Corollary~\ref{readColCont} that if monomials are ordered
by $\diag$ and the polynomials $[T]$ corresponding to straight
tableaux are ordered lexicographically by their modified column words, then
the matrix expressing the $[T]$ in terms of divided powers monomials
is in echelon form with $\pm1$'s as pivots.
\end{proof}

\begin{corollary}\label{modincr}
Suppose $p=\sum_i \alpha_i T_i$ is a linear combination of 
row-standard tableaux such that 
$p = \sum_j \beta_j [S_j]$  
where the $S_j$ are distinct straight tableaux 
and where all tableaux have the same
row-convex shape $D$.
The smallest modified column word of a tableaux in the
$S_j$'s is weakly larger (lexicographically) than the smallest modified
column word appearing in the $T_i'$.
\end{corollary}
\begin{proof}
Let $c_l$ be the column of $F(D)$ containing $l$.
Suppose that $w_{T_{i_0}}\le w_{T_{i}}$ for all $i$ and
suppose $w_{S_{j_0}}<w_{S_{j}}$ for all $j\ne j_0$--recall
by Theorem~\ref{uniqColCont} that
distinct straight tableaux have distinct modified column words.
We want to show $w_{T_{i_0}}\le w_{S_{j_0}}$. 
Now because straight tableaux have distinct modified column
words $\prod_l ({w_{S_{j_0}}}_l|c_l)$ is the smallest monomial
occurring in $p$. That means that it must appear in
$\sum_i \alpha_i T_i$ if that expression is expanded out
to a polynomial in $\Superlp LP$. But if  $w_{T_{i}}$
is always larger than $w_{S_{j_0}}$ then no monomial
as small as $\prod_l ({w_{S_{j_0}}}_l|c_l)$ can appear
in $\sum_i \alpha_i T_i$.
\end{proof}

The next section shows that any $\sum_i \alpha_i[T_i]$
can be rewritten in the above fashion.

\section{A straightening algorithm} 
\label{sect:straighten}

We produce an explicit two-rowed straightening law
for reducing any tableau to a linear combination of straight tableau.
This algorithm, {\tt straighten-tableau} shown in Figure~\ref{straighten-tableau},
starts with a tableau $T$ and returns a formal linear combination
$\sum_i \alpha_i S_i$ of straight tableau with integer coefficients such that 
$[T]=\sum_i \alpha_i [S_i]$.
In each step, the algorithm looks for a pair of rows containing a flippable inversion.
If these exist, it applies the sub-algorithm {\tt row-straighten}
in Figure~\ref{row-straighten} to ``straighten'' these two rows via
the Grosshans-Rota-Stein syzygies of Definition~\ref{syz-def}.

\begin{figure}[hbtp]
\framebox{\vbox{
{\bf Input}: A row-convex tableau $T$.
\\
{\bf Output}: $\sum_\iota \alpha_\iota S_\iota$ such that 
$[T]=\sum_\iota \alpha_\iota [S_\iota]$ 
where each $S_\iota$ is a straight tableau and $\alpha_\iota\in\Z$.
\begin{code}
\If $T$ is straight {\bf then} output $T$.
\\
{\bf else} there exists a flippable inversion in some rows $i,j$\\
\>\parbox{4.2in}{\raggedright
     Let $\sum_\kappa \beta_\kappa\cdot 
    {\aword{\mbf{v_\kappa}} \atop \aword{\mbf{w_\kappa}}}$
     be the output of $\hbox{\bf row-straighten}\!\!
         \left(\aword{{\SC t}_i} \atop \aword {{\SC t}_j}\right)$;\\
Let $N_\kappa$ be 
 $\hspace*{8pt}(\hbox{\footnotesize \# pos. letters in $\mbf{w_\kappa}$ + \# pos. letters in ${\SC T}_j$})
\cdot
 (\hbox{\footnotesize \# pos letters in ${\SC t}_{i+1}\cdots {\SC t}_{j-1}$})$.  
}\\
\>Output $\displaystyle \sum_\kappa (-1)^{N_\kappa} \beta_\kappa\cdot 
\hbox{\bf straighten-tableau}\!\!
         \left({\renewcommand{\arraystretch}{.5}
\begin{array}{c}
\aword{{\SC t}_1}\\
\aword{}\\
\aword{{\SC t}_{i-1}}\\
\aword{v_\kappa}\\
\aword{{\SC t}_{i+1}}\\
\aword{}\\
\aword{{\SC t}_{j-1}}\\
\aword{w_\kappa}\\
\aword{{\SC t}_{j+1}}\\
\aword{}\\
\end{array}}
\right)$.  
\end{code}
\caption[Algorithm straighten-tableau]{
{\bf Algorithm straighten-tableau}. \label{straighten-tableau}
${\SC t}_i$, ${\SC t}_j$, etc. are the $i$th, $j$th, etc. rows of the tableau
$T$.
}}}\end{figure}

\begin{figure}[hbtp]
\framebox{\vbox{
{\bf Input}: A two-rowed row-convex tableau 
     $T={v_{m_1} v_{m_1+1}\cdots v_{\lambda_1} \atop 
      w_{m_2} w_{m_2+1}\cdots w_{\lambda_2}}$
which is row-standard but not straight.
\\
{\bf Output}: $\sum_\kappa \alpha_\kappa\cdot  
                 {T_\kappa}$ 
such that 

Claim 1:
$[T]=
    \sum_\kappa \alpha_\kappa 
           \left[ T_\kappa \right]$ 
where $\alpha_\kappa\in\Z$ and 

Claim 2: the column word of ${\aword{v_\kappa}\atop \aword{w_\kappa}}$
is lexicographically larger than the column word of $T$.
\begin{code}
\parbox{4.5in}{
Let $c_2$ be the index of the column containing the leftmost flippable
inversion.\\
 Let $c_1$ be the smallest column such that $c_1\ge m_2$ and
either $v_{c_1-1}\lp w_{c_1}$ or $c_1-1<m_1$ (i.e. $v_{c_1-1}$ does not 
exist.)\\
 Let $c_3$ be the rightmost column such that $w_{c_2}=w_{c_3}$.}\\
\If $c_1<c_2$ {\bf then}\\
\> Let $\sum_{\iota\in I} \beta_\iota T_\iota =
       Syz_{c_2,c_2+1,\cdots,\lambda_1;c_1,c_1+1,\cdots,c_3}(T)$.\\
\>  $Expansion\gets 0$.\\
\> \For ${\iota\in I}$\\
\>\> \If $w_{T_\iota}>w_T$ {\bf then} 
              $Expansion\gets Expansion+\beta_\iota T_\iota$.\\
\>\>\> {\bf else}  
   $Expansion\gets Expansion+\beta_\iota  \cdot
           \hbox{\bf row-straighten}(T_\iota)$.\\
\> Output $Expansion$.\\
{\bf else} \Comment Comment: $c_1=c_2$. \\
\> Let $c_0$ be the leftmost column such that 
             $v_{c_0}\gp w_{c_2}$. 
            \\ \Comment Comment:  $c_0<c_2  \implies v_{c_0}= w_{c_2}$.\\
\> Let $\sum_{\iota\in I} \beta_\iota T_\iota =
       Syz_{c_0,c_0+1,\cdots,\lambda_1;c_1,c_1+1,\cdots,c_3}(T)$.\\
\>  $Expansion\gets 0$.\\
\> \For ${\iota\in I}$\\
\>\> \If $w_{T_\iota}>w_T$ {\bf then} 
              $Expansion\gets Expansion+\beta_\iota T_\iota$.\\
\>\>\> {\bf else}  
   $Expansion\gets Expansion+\beta_\iota \cdot
          \hbox{\bf row-straighten}(T_\iota)$.\\
\> Output $Expansion$.\\
\end{code}
\caption[Algorithm row-straighten]{
{\bf Algorithm row-straighten}.
If $\C L=\C L^-$, then we will always have $w_{T_\iota}>w_T$ so the
algorithm will never recurse and instead could have directly output
the expressions $Syz(T)$. The expression $Syz(T)$ is defined
in Definition~\ref{syz-def}.
 \label{row-straighten}}
}}\end{figure}

We provide an example of the straightening law below. 

\nonote{when does $\mbf w$ mean a word and when a row of a tableau?
the indexing is different.}
\begin{example}
Let $\C L = \C L^- = \{1,2,\ldots,8\}$.
In each step we shall look for a non--straight tableaux $T$
and locate two rows (say $r_1$ above $r_2$) in $T$ containing a flippable
inversion. In this example we will mark by a~$\star$
every cell in row~$r_1$ weakly right of the left most flippable
inversion in those rows and every cell in row~$r_2$ that is weakly
left of this flippable inversion and weakly right of a cutoff column $c_1$.
The cutoff $c_1$ indexes the leftmost column of row~$r_2$ such that
either $T_{r_1,c_1-1}$ does not exist or $T_{r_1,c_1-1}\lp T_{r_2,c_1}$. 
In this example, $c_1$ happens to always index the leftmost column in row~$r_2$.
We mark the remaining elements in row~$r_1$ by~$\bullet$'s.

The Grosshans-Rota-Stein syzygies (proved for the commutative case $\C L = \C L^-$
in~\cite{DRS}) says that anti-symmetrizing all the $\star$'d elements in $T$, is
the same (up to sign) as collecting all the $\star$'d elements into the 
row~$r_1$,  replacing those removed from row~$r_2$ with these $\bullet$'d
elements, and  anti-symmetrizing the $\bullet$'d elements. We shall
repeatedly apply this identity.

So, observing that the entries $4$ and $2$ form a flippable inversion, we first have,
{
\[\tblx w{\smush{\phantom{\row{5^{\bullet}}}}
\row{ ,  ,  4^\smush{\star}, 5^\smush{\star}}\\
\row{1, 3,  5,   7}\\
\row{ ,  ,  2^\smush{\star}   }\\
\row{ , 3,  8   }\\
}
=
\tblx w{\smush{\phantom{\row{5^{\bullet}}}}
\row{ ,  ,  2^\smush{\star}, 5^\smush{\star}}\\
\row{1, 3,  5,   7}\\
\row{ ,  ,  4^\smush{\star}   }\\
\row{ , 3,  8   }\\
}
-\tblx w{\smush{\phantom{\row{5^{\bullet}}}}
\row{ ,  ,  2^\smush{\star}, 4^\smush{\star}}\\
\row{1, 3,  5,   7}\\
\row{ ,  ,  5^\smush{\star}   }\\
\row{ , 3,  8   }\\
}.
\]}

\noindent
But the cells in column~3 and rows~2 and~3 of the first tableau on the
right hand side now contain a flippable inversion. We straighten as follows, 
\[
\tblx w{\smush{\phantom{\row{5^{\bullet}}}}
\row{   ,    ,  2,   5  }\\
\row{1^\smush{\bullet}, 3^\smush{\bullet},  5^\smush{\star}, 7^\smush{\star}}\\
\row{   ,    ,  4^\smush{\star}     }\\
\row{   ,   3,  8       }\\
}
=
\tblx w{\smush{\phantom{\row{5^{\bullet}}}}
\row{   ,    ,  2,   5  }\\
\row{1^\smush{\bullet}, 3^\smush{\bullet},  4^\smush{\star}, 7^\smush{\star}}\\
\row{   ,    ,  5^\smush{\star}     }\\
\row{   ,   3,  8       }\\
}
-\tblx w{\smush{\phantom{\row{5^{\bullet}}}}
\row{   ,    ,  2,   5  }\\
\row{1^\smush{\bullet}, 3^\smush{\bullet},  4^\smush{\star}, 5^\smush{\star}}\\
\row{   ,    ,  7^\smush{\star}     }\\
\row{   ,   3,  8       }\\
}
+\tblx w{\smush{\phantom{\row{5^{\bullet}}}}
\row{   ,    ,  2,   5  }\\
\row{1^\smush{\bullet}, 4^\smush{\star},  5^\smush{\star}, 7^\smush{\star}}\\
\row{   ,    ,  3^\smush{\bullet}     }\\
\row{   ,   3,  8       }\\
}
-\tblx w{\smush{\phantom{\row{5^{\bullet}}}}
\row{   ,    ,  2,   5  }\\
\row{3^\smush{\bullet}, 4^\smush{\star},  5^\smush{\star}, 7^\smush{\star}}\\
\row{   ,    ,  1^\smush{\bullet}     }\\
\row{   ,   3,  8       }\\
}{ .}
\]

\noindent
Now the first two tableaux above are straight, but the last two are not. We
straighten the next to last tableau by,
\[
\tblx w{\smush{\phantom{\row{5^{\bullet}}}}
\row{   ,    ,  2,   5  }\\
\row{1^\smush{\bullet}, 4^\smush{\star},  5^\smush{\star}, 7^\smush{\star}}\\
\row{   ,    ,  3       }\\
\row{   , 3^\smush{\star},  8       }\\
}
=
\tblx w{\smush{\phantom{\row{5^{\bullet}}}}
\row{   ,    ,  2,   5  }\\
\row{1^\smush{\bullet}, 3^\smush{\star},  5^\smush{\star}, 7^\smush{\star}}\\
\row{   ,    ,  3       }\\
\row{   , 4^\smush{\star},  8       }\\
}
-
\tblx w{\smush{\phantom{\row{5^{\bullet}}}}
\row{   ,    ,  2,   5  }\\
\row{1^\smush{\bullet}, 3^\smush{\star},  4^\smush{\star}, 7^\smush{\star}}\\
\row{   ,    ,  3       }\\
\row{   , 5^\smush{\star},  8       }\\
}
+
\tblx w{\smush{\phantom{\row{5^{\bullet}}}}
\row{   ,    ,  2,   5  }\\
\row{1^\smush{\bullet}, 3^\smush{\star},  4^\smush{\star}, 5^\smush{\star}}\\
\row{   ,    ,  3       }\\
\row{   , 7^\smush{\star},  8       }\\
}
+\tblx w{\smush{\phantom{\row{5^{\bullet}}}}
\row{   ,    ,  2,   5  }\\
\row{3^\smush{\star}, 4^\smush{\star},  5^\smush{\star}, 7^\smush{\star}}\\
\row{   ,    ,  3       }\\
\row{   , 1^\smush{\bullet},  8       }\\
}
\]

\noindent and the last tableau by
\[
-\tblx w{\smush{\phantom{\row{5^{\bullet}}}}
\row{   ,    ,  2^\smush{\star}, 5^\smush{\star}  }\\
\row{3  , 4  ,  5  , 7  }\\
\row{   ,    ,  1^\smush{\star}     }\\
\row{   ,   3,  8       }\\
}
=
-\tblx w{\smush{\phantom{\row{5^{\bullet}}}}
\row{   ,    ,  1^\smush{\star}, 5^\smush{\star}  }\\
\row{3  , 4  ,  5  , 7  }\\
\row{   ,    ,  2^\smush{\star}     }\\
\row{   ,   3,  8       }\\
}
+\tblx w{\smush{\phantom{\row{5^{\bullet}}}}
\row{   ,    ,  1^\smush{\star}, 2^\smush{\star}  }\\
\row{3  , 4  ,  5  , 7  }\\
\row{   ,    ,  5^\smush{\star}     }\\
\row{   ,   3,  8       }\\
}{ ,}
\]
so

{\footnotesize
\[\displaylines{\tblx w{\smush{\phantom{\row{5^{}}}}
\row{ ,  ,  4^\smush{}, 5^\smush{}}\\
\row{1, 3,  5,   7}\\
\row{ ,  ,  2^\smush{}   }\\
\row{ , 3,  8   }\\
}
=
-\tblx w{\smush{\phantom{\row{5^{}}}}
\row{ ,  ,  2^\smush{}, 4^\smush{}}\\
\row{1, 3,  5,   7}\\
\row{ ,  ,  5^\smush{}   }\\
\row{ , 3,  8   }\\
}
+\tblx w{\smush{\phantom{\row{5^{}}}}
\row{   ,    ,  2,   5  }\\
\row{1^\smush{}, 3^\smush{},  4^\smush{}, 7^\smush{}}\\
\row{   ,    ,  5^\smush{}     }\\
\row{   ,   3,  8       }\\
}
-\tblx w{\smush{\phantom{\row{5^{}}}}
\row{   ,    ,  2,   5  }\\
\row{1^\smush{}, 3^\smush{},  4^\smush{}, 5^\smush{}}\\
\row{   ,    ,  7^\smush{}     }\\
\row{   ,   3,  8       }\\
}
+\tblx w{\smush{\phantom{\row{5^{}}}}
\row{   ,    ,  2,   5  }\\
\row{1^\smush{}, 3^\smush{},  5^\smush{}, 7^\smush{}}\\
\row{   ,    ,  3       }\\
\row{   , 4^\smush{},  8       }\\
}+\cr
\hfill 
-\tblx w{\smush{\phantom{\row{5^{}}}}
\row{   ,    ,  2,   5  }\\
\row{1^\smush{}, 3^\smush{},  4^\smush{}, 7^\smush{}}\\
\row{   ,    ,  3       }\\
\row{   , 5^\smush{},  8       }\\
}+
\tblx w{\smush{\phantom{\row{5^{}}}}
\row{   ,    ,  2,   5  }\\
\row{1^\smush{}, 3^\smush{},  4^\smush{}, 5^\smush{}}\\
\row{   ,    ,  3       }\\
\row{   , 7^\smush{},  8       }\\
}
+\tblx w{\smush{\phantom{\row{5^{}}}}
\row{   ,    ,  2,   5  }\\
\row{3^\smush{}, 4^\smush{},  5^\smush{}, 7^\smush{}}\\
\row{   ,    ,  3       }\\
\row{   , 1^\smush{},  8       }\\
}
-\tblx w{\smush{\phantom{\row{5^{}}}}
\row{   ,    ,  1^\smush{}, 5^\smush{}  }\\
\row{3  , 4  ,  5  , 7  }\\
\row{   ,    ,  2^\smush{}     }\\
\row{   ,   3,  8       }\\
}
+\tblx w{\smush{\phantom{\row{5^{}}}}
\row{   ,    ,  1^\smush{}, 2^\smush{}  }\\
\row{3  , 4  ,  5  , 7  }\\
\row{   ,    ,  5^\smush{}     }\\
\row{   ,   3,  8       }\\
}.\qquad\cr
}\]
}
\end{example}

The first step in verifying  Algorithm~{\tt straighten-tableau} is
to prove that when $T$ is replaced with $\sum_i \beta_i T_i$ 
by Algorithm~{\tt row-straighten} we have
$[T]=\sum_i \beta_i [T_i]$. The second step involves showing that 
each $T_i$ is somehow closer to being straight than was $T$. 
The first of these facts is an immediate consequence of the
correctness of Algorithm~{\tt row-straighten}. 
This will come down to
 verifying the identities used in the preceding example.
The second will follow
from the correctness of {\tt row-straighten} and the fact, 
proved in Proposition~\ref{lowerColWord}, that
given a tableau $T$ and another tableau $T'$ differing only in
two rows $i,j$, then the column word of the two-rowed subtableaux
consisting of rows $i,j$ of $T$ is less than the corresponding column
word determined by $T'$ iff $c_T<c_{T'}$.

The proof of  Algorithm~{\tt straighten-tableau} thus depends solely
on the correctness of Algorithm~{\tt row-straighten}.
We will prove both claimed properties of Algorithm~{\tt row-straighten} for 
each of the two cases appearing in the algorithm.
First we will produce the ``determinantal'' identities
that will be used in Algorithm~{\tt row-straighten}.

Define a \defi{shuffle} of a word $\mbf w=w_1,\ldots,w_n$ into parts of length
$k,k'=n-k$ to be an ordered pair of words $\mbf{w'}$  and $\mbf{w''}$  of 
$\mbf w$ having lengths $k$ and $k'$ respectively, such that 
$\mbf{w'}$  and $\mbf{w''}$
can be found as a pair of disjoint subwords of $\mbf w$. Neither 
$\mbf{w'}$  nor $\mbf{w''}$ need be contiguous as a subword of $\mbf w$.
 When $\mbf w=1,\ldots, n$ 
a shuffle amounts to a permutation $\sigma$ of the 
index set $1,\ldots,n$ such that $\sigma_1<\cdots<\sigma_k$ and 
$\sigma_{k+1}<\cdots<\sigma_{k+k'}$. Generalizing the length of
a permutation we define the \defi{shuffle signature},
$\shflsign(\mbf\omega_1,\ldots,\omega_k)$, of a word $\mbf\omega$ to be the
number of pairs $1\le i<j\le k$ such that $\omega_i>\omega_j$ and  
$|\omega_i|=|\omega_j|=1$.


\begin{definition}\label{syz-def}
{\def\s{\sigma}
Let $i,j,k,l$ be nonnegative integers.
Fix a two-rowed row-convex shape $D$ by
specifying the starting and ending columns, $1$ through $i+j+l$ and 
$m$ through $m+l+k-1$ of the top and bottom rows respectively. 
For convenience we have let the leftmost column index of the top row be~$1$,
but the column indices can of course be shifted left or right by any integer.
With the above convention, we could have $m\le 0$, this produces a skew shape.

Let $T={\aword{\mbf v} \atop \aword{\mbf w}}$ be a two-rowed row-convex tableau of shape $D$.
Fix two sequences of column indices,
$c_1<\cdots <c_j$ starting at column $i+l+1$ and ending at $i+l+j$
and 
$m\le c'_1<\cdots<c'_l\le k+l$ starting at column $c_1$ and ending
at $c_1+l-1$.

{\em Define} $Syz_{c_1,\ldots,c_j;\,\,c'_1,\ldots,c'_l}(T)$ to be the formal linear combination
\begin{equation}\label{SyzExpr}
\sum_{{  {\scriptstyle \rm all\ nontrivial\ shuffles\ } 
\atop 
  {{\scriptstyle
     y_{\sigma(1)}  \ldots y_{\sigma(l)} ;\  y_{\sigma(l+1)} \ldots y_{\sigma(l+j)}    }
\atop
    {\scriptstyle {\rm of}\ 
     w_{\sigma(c'_1)} \ldots w_{\sigma(c'_l)}   v_{\sigma(c_1)}  \ldots v_{\sigma(c_j)}  }
  }} }
\!\!\!\!\!\!  -{\alpha_\sigma\over \alpha_e} \cdot T'_{\sigma}\,
+
\sum_{{  {\scriptstyle \rm all\ shuffles\ } 
\atop 
  {{\scriptstyle
     x_{\tau(1)} \ldots x_{\tau(i)} ;\  x_{\tau({i+1})}  \ldots x_{\tau({i+l})} }
\atop
    {\scriptstyle {\rm of}\  v_1 v_2 \cdots v_{c_1-1}} 
}}}
\!\!\!\!\!\!  {\beta_\tau \over \alpha_e} \cdot T''_{\tau},
 \end{equation}
\def\s{\sigma}
where $T'_{\sigma}$, $T''_{\tau}$, $\alpha_\sigma$ and $\beta_\tau$ are defined as follows:

First, define words 
\[\begin{array}{cc}
 \mbf x = v_1,\ldots,v_{c_1-1}\,, &
            u_1,\ldots,u_{j+l} = w_{c'_1},\ldots,w_{c'_l},v_{c_1},\ldots,v_{c_j}\,, \\  
 \mbf{x'_\tau}   =  v_{\tau(1)},\ldots,v_{\tau(i)}    \,, &
 \mbf{u'_\sigma} =  u_{\s(1)},\ldots,u_{\s(l)}    \,, \\
 \mbf{x''_\tau}  =  v_{\tau(i+1)},\ldots,v_{\tau(i+l)}   \,,  &
 \mbf{u''_\sigma}=  u_{\s(l+1)},\ldots,u_{\s(l+j)}   \,, \\
\mbf{z_1} = w_m,\ldots,w_{c'_1-1} & \mbf{z_2} = w_{c'_l+1},\ldots,w_{k+l-1}\, \\ 
{\rm\ and} & \mbf z = \mbf{z_1},\mbf{z_2}\,. 
\end{array}\]

Define the tableau $T'_{\sigma}$ to be the tableau obtained by sorting the rows of
\[\rsktblx n{
\srow{\cdots\aword{\mbf{x}}\cdots,.... \cdots\cdots\aword{\mbf{u''_\sigma}}\cdots\cdots}\\ 
\srow{,..  \aword{\mbf{z_1}},. \aword{\mbf{u'_\sigma}},. \aword{\mbf{z_2}}  }\\ 
}\] 
and $T''_{\tau}$ to be to be the result of sorting the rows of\nonote{in revision switch $u$ and $x'_\tau$
here and in sorting sign for $N_\tau$ and add the $|u|\cdot|x'_\tau|$ required for switching to $N_\tau$}
\[\rsktblx n{
\srow{\cdots\cdots\cdots\aword{\mbf{u}}\cdots\cdots\cdots,..  \aword{\mbf{x'_\tau}}\cdots }\\
\srow{     ,...        \cdots\aword{\mbf{x''_\tau}},... \cdots\aword{\mbf z}\cdots  }\\ 
}.\]

We define  
\[ \alpha_\sigma =
  (-1)^{N_1(\tau)}\cdot{\c(\mbf x, \mbf{u''_\sigma})! 
                          \over  
      \c(\mbf x)! \cdot \c(\mbf{u''_\sigma})! }
\cdot 
{\c(\mbf{z_1},u'_\sigma,\mbf{z_2})! 
                         \over 
\c(\mbf{z_1},\mbf{z_2})! \cdot \c(\mbf{u'_\sigma})!}
\]
with 
$N_1(\tau)=|\mbf{z_1}| \cdot |\mbf{u'_\sigma}| + 
|\mbf{z}|\cdot (i+j+l) +
|u'_\sigma|\cdot(i+j+l) +
\shflsign(\mbf{u''_\sigma},\mbf{u'_\sigma}) +
m(\mbf{x}, \mbf{u''_\sigma}) + 
m(\mbf{z_1}, \mbf{u'_\sigma}, \mbf{z_2}  )
$
where $m(\mbf\omega)=\shflsign(\mbf\omega)+{\hbox{\# neg. letters in $\mbf\omega$}\choose 2}$
and we define 
\[ \beta_\tau =
(-1)^{l}\cdot(-1)^{N_2(\tau)}\cdot
{\c(\mbf u, \mbf{x'_\tau})! 
\over  
      \c(\mbf u)! \cdot \c(\mbf{x'_\tau})! }
\cdot
{\c(\mbf{x''_\tau}, \mbf{z})!
\over 
\c(\mbf{x''_\tau})! \cdot \c(\mbf{z})!}
\]
with
$N_2(\tau)= |\mbf{z}|\cdot (i+j+l)  +
|\mbf{x''_\tau}|\cdot \big( (i+j+l) - |\mbf{u}|\big) +
\shflsign(\mbf{x'_\tau},\mbf{x''_\tau}) +
m(\mbf{u}, \mbf{x'_\tau} ) +
m(\mbf{x''_\tau}, \mbf z) 
$.

Call the cells in columns $c_1,\ldots,c_j$ of the top row and 
the cells in columns $c'_1,\ldots,c'_l$ of the bottom row of $T$ \defi{marked cells}.
If no positive letter appears in both marked and unmarked cells in the top row of $T$
and no positive letter appears in both marked and unmarked cells in the bottom row of $T$
then $\alpha_e=\pm1$ so the identity holds over $\Z$. 
}\end{definition}


\begin{proposition}\label{polzn-ident}
Let $T$ be a tableau of two-rowed, row-convex shape $D$ whose top row contains its bottom row
Without loss of generality assume the top row has leftmost column index 1.
If column indices $c_1,\ldots,c_j$ and $c'_1,\ldots, c'_l$ are chosen as in
Definition~\ref{syz-def}, then
\begin{equation}\label{expansionidentity}
 [T] =  \sum_S a_S\cdot[S] 
\end{equation}
where $\sum_S a_S\cdot S = Syz_{c_1,\ldots,c_j;\,\,c'_1,\ldots,c'_l}(T)$
and both sums are over all (w.l.o.g. row-standard) tableau of shape $D$.

The identity~\ref{expansionidentity} also holds when $D$ is a two-rowed skew shape 
(i.e. $M\le0$ in Definition~\ref{syz-def}) and
$l+j$ exceeds the number of columns in $D$. In this case, the right-hand summand in 
expression~\ref{SyzExpr} vanishes.
\end{proposition}

This result follows from the more general Theorem~10 of~\cite{GRS}
\nonote{in br-rule sectn, you could produce the ident's
here as a cyclic module over the univ env. alg.}
but for completeness, we sketch a proof relying on positive letters and the polarization operators of 
Section~\ref{sect-br-rule} that bypasses 
the Hopf algebra techniques of~\cite{GRS}. The proof provides a much
simpler though less explicit definition of the expressions $Syz_{\mbf c, \mbf c'}$.

\begin{proof}
First, we prove the following proposition directly by checking that the monomials
arising from the expansion of each expression have the same coefficients. 
Details may be found in~\cite{thesis}.
\begin{proposition}\label{pos-ident}
Let $a,b,c$ be positive letters.  Let $i,j,k,l$ be nonnegative integers.
Let $1,2,\ldots, i+j+l$ be negative letters.
Fix a two-rowed row-convex shape $D$ whose top row contains its bottom row by
specifying the starting and ending columns, $1$ through $i+j+1$ and 
$m$ through $m+l+k-1$ of the top and bottom rows respectively. The following
identity holds for tableaux of shape $D$:
\refstepcounter{equation}\label{posversion}
\begin{displaymath}
\displaylines{
\qquad
\bitblx n{
\srow{a^{(i+l)}, b^{(j)}} &\vbar \srow{1, 2,       ,   \ldots, \ldots,      \ldots,  \sst i+j+l}\\
\srow{b^{(l)}, c^{(k)}}   &\vbar \srow{  ,      m, \sst m+1, \ldots, \sst m+l+k-1}\\
}\hfill\cr
\hfill=(-1)^l
\bitblx n{
\srow{b^{(j+l)}, a^{(i)}} &\vbar \srow{1, 2,       ,   \ldots, \ldots,      \ldots,  \sst i+j+l}\\
\srow{a^{(l)}, c^{(k)}}   &\vbar \srow{   ,      m, \sst m+1, \ldots, \sst m+l+k-1}\\
}\qquad(\theequation)}
\end{displaymath}
 where 
$\left( {\aword{{\SC s}_1}\atop\aword{{\SC s}_2}} \bigm| {\aword{{\SC t}_1}\atop\aword{{\SC t}_2}} \right)
=
( {\aword{{\SC s}_1} \mid \aword{{\SC t}_1} })
( {\aword{{\SC s}_2} \mid \aword{{\SC t}_2} })$
and where $l^{(k)}=l^k/k!$.
\end{proposition}

Choosing $\mbf v= v_1,\ldots,v_{i+l}$, $\mbf w = w_1,\ldots, w_k$, and 
$\mbf w = u_1,\ldots, u_{j+l}$.
and applying the product,
\[{1\over  \c(\mbf v)!\c(\mbf w)!\c(\mbf u)! }
D_{v_1,a} D_{v_2,a} \cdots D_{v_{i+l},a} \ 
D_{w_1,c}  \cdots D_{w_{k},c}          \ 
D_{u_1,b}  \cdots D_{u_{j+l},b} \, ,\]
of polarizations to the identity in Proposition~\ref{pos-ident}
completes the proof when $m$ in Definition~\ref{syz-def} is positive.
For $m\le0$, we recover the identity used traditionally to straighten 
skew tableaux. A proof follows by recognizing that the right-hand side
of equation~\ref{posversion} vanishes when $j+l>i+j+l-m+1$.
\end{proof}

We have just verified that any formal linear combination of tableau with
integer coefficients produced by
Algorithm~{\tt row-straighten} satisfies Claim~1 made in
the algorithm specifications; we now go to work on 
the heart of the proof, namely Claim~2.
In the process of proving Claim~2, we will check that 
{\tt row-straighten} terminates.

\begin{proposition}\label{lowerColWord}
Let $D$ be the row-convex shape $(\lambda_1,\lambda_2)/(m_1,m_2)$ with $m_1<m_2$.

Given a non-straight, row standard, two rowed, row-convex tableau,  $T$,
of shape $D$
Algorithm~{\tt row-straighten} produces a formal linear
combination of tableaux each of which has a lexicographically
larger column word than $c_T$.
\end{proposition}

\begin{proof}
The proof is by induction on $c_T$.

Suppose that $T$ is the tableau \[
 \rsktblx n{
\row{\llap{Column \#\ \ } 
   \sst  m_1, \sst m_1+1, \sst m_2, \sst m_{2+1} , , \sst \lambda_2, \sst \lambda_1}
\row{\vrule height0pt depth 20pt width0pt}\\
\row{v_{m_1}, v_{m_1+1},        ,          , \ldots,              , v_{\lambda_1}  } \\
\row{       ,          , w_{m_2}, w_{m_2+1}, \ldots, w_{\lambda_2}  }   \\
}.\] 

We set up a straightening syzygy that expresses $T$ in terms of tableau $T_i$
such that $c_T$ is always lexically smaller than $c_{T_i}$. 
We first describe the structure of $T$ with respect to its leftmost flippable
inversion. Let $q$ be minimal such that if column~$q+1$
were to contain an inversion, then that inversion would be flippable.
Thus $q= \min_{m_2\le i\le \lambda_2\atop v_{i-1} \lp w_i} i-1 $.
The presence of a flippable inversion guarantees that $q$ exists.
The value of $c_1$ in Algorithm~{\tt row-straighten} is $q+1$.
The first column strictly right of $c_1-1$ that actually {\em has} an inversion
has index
$  c_2 = \min_{c_1\le j\le\lambda_2\atop v_j\lp w_j} j $.
This inversion must be flippable, thus $c_2$ indexes the column of
the leftmost flippable inversion in the tableau.

Two cases arise in the algorithm namely $c_1<c_2$ and $c_1=c_2$. The  pictures 
in Figures~\ref{setup1} and~\ref{setup2} outline
these situations. The symbol, ``$\bullet$'', indicates a cell in the
diagram. An arrow from one cell to another indicates that the contents
of the first cell are larger than the contents of the second. 
The decoration of an arrow by ``$-$'' (respectively ``$+$'') indicates
that the contents of the cells at either end are may be equal
if these contents are negatively (respectively positively) signed.
Sequences of cells surrounded by parentheses or braces may be omitted.

\begin{figure}[tbp]
\framebox{\vbox{
\[\rsktblx n{
\smush{\row{\phantom{aaa..}}}
\smash{\smush{\row{
     \up{m_1},       ,      ,       ,       ,       ,       , \up{c_1-1},    \up{\ c_1}, ,       , \up{c_2-1}, \up{\ c_2},       \up{\lambda_2},           \up{\lambda_1}  }}}
\smash{\smush{\row{
             ,       ,      ,    \smush{\overbrace{\row{ , ,  }}},   }}}
\row{     \bu, \ldots,   \bu,    \bu,    \bu, \ldots,    \bu,   \bu,      \bu,   \bu, \ldots,   \bu,        \bu,  \llap{\ldots\,}\bu,            \llap{\ldots\,}\bu}\\
\smush{\row{\lprn,   , \rprn,       ,       ,       ,       ,      ,         , \lprn,       ,             \rprn,   }}
\smash{\smush{\row{ 
      ,       ,      ,      , \bseam, \bseam, \kern-25pt\ldots, \bseam, \bnwa\buam, \buam, \ldots, \buam,     \bdap, }}}\\
\smash{\smush{\row{ 
      ,       ,      ,      ,  \smush{\underbrace{\row{ , ,  }}}, }}}
\row{         ,      ,      ,       ,    \bu,    \bu, \ldots,   \bu,    \bu,   \bu, \ldots,   \bu,        \bu,     \llap{\ldots\,}\bu} \\
\smash{\smush{\row{
      ,       ,      ,      , \down{m_2},   ,       ,       ,      ,      ,       ,      ,           ,  \llap{\down{c_3}\ \ \ }       }}}\\
}
\]
\caption[Setup for two-row straightening, Case~I]{
{\bf Case I}: $c_1<c_2$: Relations between entries in a two-row tableau being straightened by
Algorithm~{\tt row-straighten}.

All entries in the bottom row from $c_2$ through $c_3$ are 
equal but distinct
from any entry in column $c_3+1$ of that row.  
 \label{setup1} }
}}\end{figure}

If no positive letter appears multiple times in $T$, then 
Cases~I and~II can be treated simultaneously.
We begin with Case~I.

Apply Corollary~\ref{polzn-ident} to write
\begin{equation}\label{spanning:intoAB}
\tblx n{
\smash{\smush{
\row{       ,        ,       ,        ,       , \smush{\kern-12pt\overline{\row{ , , , , }}}  } }}
\row{v_{m_1},        , \ldots,  \ldots, \ldots, v_{c_2},       , \ldots,  , v_{\lambda_1}  } \\
\smash{\smush{
\row{       ,        ,       , \smush{\kern-12pt\underline{\row{ , , , , }}}      } }}
\row{       , w_{m_2}, \ldots, w_{c_1}, \ldots, w_{c_2}, \ldots, w_{c_3}, \ldots } \\
} = \mathop{Syz}\limits_{c_2,\ldots,\lambda_2;\, c_1\ldots c_3}  =  B + A  
\end{equation}
where $B$ (respectively $A$) is the first (respectively second) summation in the 
$Syz_{c_2,\ldots,\lambda_2;\, c_1\ldots c_3}(T)$ as defined in expression~\ref{SyzExpr}.
The over/underlines are visual aids which indicate the ``marked'' entries used to
define $Syz(T)$.

It suffices to show that each tableau appearing in $A$ or $B$ has lexically larger
column word than $c_T$.

Suppose $T'$ appears in $A$. We can write
\[ T'=\rsktblx n{
\smush{\row{ \kern-20pt\raisebox{-4.75pt}{\makebox[0pt][l]{\frame{\row{ , , , ,..... }}}}}}
\row{x_{m_1},   \ldots,  x_{t-1}, w_{c_1}, \ldots, w_{c_3},   v_{c_2}, \ldots, \ldots, \ldots,.... v_{\lambda_1} }\\
\smush{\row{ , \kern-20pt\raisebox{-4.75pt}{\makebox[0pt][l]{\frame{\row{ , ,.... , ,. ,  }}}}}}
\row{ ,
w_{m_2}, \ldots,....      w_{c_1-1},    y_1,  \ldots,.... y_{c_2-t}, w_{c_3+1}, \ldots, w_{\lambda_2} }\\
}, \]
where $x_{m_1}\ldots x_{t-1}  ;\, y_1\ldots y_{c_2-t}$  is a shuffle 
of $v_{m_1}\ldots v_r$, where $t=c_2+c_1-c_3-1$, and 
where the boxed entries must be sorted in order to give a row-standard tableau.
Denote the entries in the bottom row by $z_{m_2},\ldots,z_{\lambda_2}$
and the entries in the top row by $\varpi_{m_1},\ldots,\varpi_{\lambda_1}$.

To check that only the boxed entries need to be sorted in the to row, 
it suffices to observe that the $x_i$ are taken from $v_{m_1},\ldots,v_r$
and that $v_{c_2}\gp w_{c_2} = w_{c_3}$.
Checking the bottom row, it suffices to note that $w_{c_3+1}\gp v_{c_2-1}$.

Now let $k+1$ index the leftmost column in the top row in which $T'$ differs from $T$.
In fact, $k=\min_{m_1\le i\le t \atop \varpi_i\ne v_i} i-1$
which follows from being in Case~I:                                                    
Suppose $\varpi_i=v_i$ for $m_1\le i<t=c_1-1-c_3+c_2$.
Since by construction $v_{t-1}\lp v_{c_1-1} < w_{c_1}$, we
have that  $x_i=v_i$ for all $i$ as above and thus  the boxed elements
in the top row are already in order. But then $\varpi_t=w_{c_1}\ne v_t$,
since by Case~I $v_{c_1-1}<w_{c_1}$. So $k<t$.

Now we examine the column words.
Our construction shows $\mbf\varpi \ge \mbf v$, so by the preceding 
paragraph $\mbf\varpi>\mbf v$. 
So if $k+1<m_2$ we
conclude directly that $c_{T'}$ is lexically larger than $c_T$.

Suppose that $k+1\ge m_2$.
We show that
$v_{k+1}\le y_1$. Suppose to the contrary that $v_{k+1}> y_1$. Since $y_1$ comes
from $v_{m_1},\ldots,v_{r}$ this says that $y_1=v_j$ for some $j\le k$
and $y_1\ne v_{j'}$ for $j'>k$. Now the upper row of $T'$ still contains 
$v_1\ldots v_k$ even though a $y_1$ has been removed to the bottom row.
But this implies that $y_1$ also appears in $w_{c_1}\ldots w_{c_2}$ which
is impossible since $w_{c_1}\gp v_{c_1-1}\gp v_{k+1}>y_1$. 

Thus since $k+1\le t<c_1$ the diagram for Case~I 
shows that $w_{k+1}< v_{k+1}$, hence $w_{k+1}<y_1$. So, after sorting,
we find that
$z_{m_2}=w_{m_2},\,z_{m_2+1}=w_{m_2+1},\, \ldots,z_{k+1}=w_{k+1}$.
So in tableaux $T$ and $T'$, the columns $m_1,\ldots,k$ agree as
does the bottom entry of column $k+1$. But the top entry in
column $k+1$ is larger in $T'$ than in $T$. Hence
$c_{T'}$ is lexically larger than $c_T$. 

\medbreak

At last we deal with tableaux appearing in $B$ in equation~\ref{spanning:intoAB}.
Recall that tableaux  in $B$ arise from nontrivially shuffling 
 the over/underlined entries and then resorting the rows. 
Let $\mbf u =w_{c_1}\cdots w_{c_3} v_{c_2}\cdots v_{\lambda_1}$.
Let $\mbf u'',\mbf u'$ be a shuffle of $\mbf u$ into two parts of size
$\lambda_1-c_2+1$ and $c_3-c_1+1$ respectively.
Since $w_{c_1}\gp v_{c_1-1}$, such a tableau will look like
\[ T'=\rsktblx n{
  \smush{\row{
       ,        ,       ,    ,     
          \kern-25pt\raisebox{-4.75pt}{\makebox[0pt][l]{\frame{\row{ , , , , , ,aaaa }}}}}}
  \row{
v_{m_1},        , \ldots, v_{c_1-1}, v_{c_1},  \ldots, v_r,         \word{\aword{\mbf u''}}, , , ,  }\\
  \smush{\row{
       ,        ,       ,    , 
          \kern-25pt\raisebox{-4.75pt}{\makebox[0pt][l]{\frame{\row{ , , , , ,aaaa }}}}}}
  \row{
       , w_{m_2}, \ldots, w_{c_1-1},  \kern-8pt\word{\aword{\mbf u'}}, , ,   , w_{c_3+1}, \ldots, w_{\lambda_2} }\\
  }, \]
where, as before, the boxed elements must be sorted so that $T'$ will be row standard.
Again denote the top and bottom rows of $T'$ by 
$\varpi_{m_1},\ldots, \varpi_{\lambda_1}$
and $z_{m_2},\ldots, z_{\lambda_2}$.

Now let $k+1$ be the leftmost column in which the bottom rows of $T$ and $T'$ disagree.
We claim $c_1\le k+1 \le c_3$ and
that $z_{k+1}>w_{k+1}$.
By construction 
$\mbf u'\ne w_{c_1}\cdots w_{c_3}$. 
Write $\mbf u'=\mbf u'_{q+1}\cdots \mbf u'_{c_3}$.
Let $j$ be minimal such that $\mbf u'_j\ne w_j$. Because 
$v_{c_2}\gp w_{c_2}=w_{c_3}$, this
implies that $\mbf u'_j > w_j$. 
But since also $w_{c_3+1}>w_{c_3}$, 
we find $z_{j}>w_{j}$. and $z_i=w_i$ for all $i<j$, so $k+1=j$.

Subcase~1. Suppose $k<c_2-1$.
Any letter appearing in the multiset difference 
$\mbf u''-\{\{v_{c_2},\ldots,v_{\lambda_1}\}\}$ is $(\gp)$-greater than
$w_{k+1}$.
But since $k<c_2-1$, the picture of case~I 
shows that $v_{k+1}\lm w_{k+1}$,
this means that on resorting the boxed elements of the top row,
every element in $\mbf u'$ stays in column $k+2$ or higher.
Hence columns $m_1\ldots k$ agree in $T$ and $T'$. 
But the bottom of column $k+1$ is larger in $T'$ than $T$. Thus
$c_{T'}$ is lexically larger than $c_T$.

The above argument also generates the  fact (unused in this proof, but 
see the comment after Corollary~\ref{modincrtwo}) that, 
$T$ and $T'$ agree in the top element of column $k+1$.

Subcase~2: suppose that $k\ge c_2-1$. This says that the bottom rows of $T,T'$ agree 
at least through column $c_2-1$.
Since $v_{c_2-1}\lp w_{c_2}$, we have immediately that the top rows of $T,T'$ agree through
column $c_2-1$. Now either the bottom of column $c_2$ changes (hence increases) so
$c_{T'}$ is lexically larger than $c_T$ and we are done or the
the number of positive letters in the bottom row that equal $w_{c_2}$ decreases.  
In the latter case, not only do the tableaux $T,T'$ agree up to column $c_2-1$ but
$T'$ still has 
a flippable inversion in column $c_2$ since the entry in the top of that column
now equals the positive letter $w_{c_2}$ that remains at the bottom.

We repeat the straightening law on $T'$, producing some tableaux with
lexicographically larger modified column words and 
some tableaux that are unchanged in  columns smaller than $c_2$ and 
unchanged at the
bottom of column $c_2$ but which
have fewer copies of $w_{c_2}$ in their bottom rows.
Eventually, we must run out of positive letters equal to $w_{c_2}$ in
the bottom row and so eventually the modified column word increases.

\medbreak

We now treat Case~II. Here $c_1=c_2$. 

We replace $[T]$ with $Syz_{c_0,\ldots,\lambda_1 ;\, c_1,\ldots, c_3 }$
where $c_0$ is minimal such that $v_{c_0} \gp w_{c_2}$.
\nonote{check the other uses of $Syz$ to make sure the words are in the right order.}
As before, Corollary~\ref{polzn-ident} lets us write
\begin{equation}\label{case2:spanning:intoAB}
\tblx n{
\smash{\smush{
\row{       ,        ,       , \smush{\kern-12pt\overline{\row{ , , , , , , }}}  } }}
\smash{\smush{
\row{       ,        ,       , 
    \smush{\kern-12pt\overline{\row{ {\vrule height 11pt depth 0pt width 0pt},   }}}  } }}
\row{v_{m_1},        , \ldots,  v_{c_0}, \ldots, v_{c_2},       , \ldots,  , v_{\lambda}  } \\
\smash{\smush{
\row{       ,        ,       ,         ,       , \smush{\kern-12pt\underline{\row{  , , }}}      } }}
\smash{\smush{
\row{       ,        ,       ,         ,       , 
    \smush{\kern-12pt\underline{\row{  {\vrule height 0pt depth 6pt width 0pt}, , }}}      } }}
\row{       , w_{m_2}, \ldots, w_{c_1}, \ldots, w_{c_2}, \ldots, w_{c_3}, \ldots } \\
} =  B + A  \end{equation}
just as in equation~\ref{spanning:intoAB}
The entries that have been marked twice are 
positive letters all equal to each other.

Suppose $T'$ appears in $A$. 
To this purpose, let $w=v_{m_1}\ldots v_r$ be the word being shuffled.
It is easily verified that
\[ T'=\rsktblx n{
\row{x_{m_1},   \ldots,  x_{t-1}, v_{c_0}\rlap{$=$}, \ldots, \llap{$=$}w_{c_3},   v_{c_2}, 
      \ldots, , ,... \ldots,.... v_{\lambda_1} }\\
\smush{\row{ , \kern-20pt\raisebox{-4.75pt}{\makebox[0pt][l]{\frame{\row{ , , , , , ,. , }}}}}}
\row{ ,
w_{m_2}, \ldots,  , \ldots, w_{c_2-1},    y_1,  \ldots,.... y_{c_0-t}, w_{c_3+1}, \ldots, w_{\lambda_2} }\\
}, \]
where $t=c_0-c_3+c_2-1$,    $x_{m_1}\ldots x_{t-1} ;\, y_1\ldots y_{c_0-t}$   is a shuffle of  
$v_{m_1},\ldots,v_{c_0-1}$  and
where the boxed entries must be sorted in order to get a row-standard tableau.
Maintain the notation $T'={\aword{\mbf \varpi} \atop \aword{\mbf z}}$.

\begin{figure}[tbp]
\framebox{\vbox{
\[\rsktblx n{
\smush{\row{\phantom{aaaa}}}
\smash{\smush{\row{
     \up{m_1},       ,      ,      ,     ,       ,       , \up{c_2-1}, \up{c_2},    ,       , \up{\lambda_2},          , \up{\lambda_1}  }}}
\smash{\smush{\row{
             ,       ,      ,    \smush{\overbrace{\row{ , ,  }}},   }}}
\row{     \bu, \ldots,   \bu,    \bu, \bu, \ldots,    \bu,   \bu,        \bu, \bu, \ldots,            \bu, \ldots, \bu}\\
\smush{\row{\lprn,   , \rprn,       ,    ,       ,       ,      ,  }}
\smash{\smush{\row{ 
      ,       ,      ,       , \bseam, \bseam, \kern-25pt\ldots, \bseam, \bnwap\bdap,  }}}\\
\smash{\smush{\row{ 
      ,       ,      ,       ,  \smush{\underbrace{\row{ , ,  }}}, }}}
\row{ ,       ,      ,      ,   \bu,  \bu, \ldots,    \bu,       \bu,    \bu, \ldots, \bu,     } \\
\smash{\smush{\row{
      ,       ,      ,      , \down{m_2}, ,      ,       ,          ,       , \down{c_3},    ,          ,         }}}\\
}
\]
\caption[Setup for two-row straightening, Case~II]{
{\bf Case II}: $c_1=c_2$: Relations between entries in a two-row tableau being straightened by
Algorithm~{\tt row-straighten}. 

All entries in the bottom row from $c_2$ through $c_3$ are equal but distinct
from any entry in column $c_3+1$ of that row.  If $c_0<c_2$, then the entries
in the top row that equal the bottom row entry in column~$c_2$ must start at
$c_0$ and extend at least as far as $c_2-1$.
 \label{setup2} }
}}\end{figure}

Since $t<c_0$ we have $v_{c_0}\ne v_t$ and thus
$ k=\min_{m_1\le i\le t \atop \varpi_i\ne v_i} i-1 $
is well defined.
As in Case~I,  if $k+1<m_2$ we
conclude directly that $c_{T'}$ is lexically larger than $c_T$.

Suppose that $k+1\ge m_2$. 
We show that
$v_{k+1}\le y_1$. Suppose to the contrary that $v_{k+1}> y_1$. Since $y_1$ comes
from $v_{m_1},\ldots,v_{c_0-1}$ this says that $y_1=v_j$ for some $j\le k$
and $y_1\ne v_{j'}$ for $j'>k$. But this says that if the letter $y_1$ 
occurs in the 
$y_1\cdots y_{c_0-t}$ part of the shuffle, then the $x_{m_1}\cdots x_k$ part
cannot start with $v_{m_1}\cdots v_k$--contradiction.

Thus since the diagram for Case~II 
shows that $w_{k+1}< v_{k+1}$, we find
 $w_{k+1}<y_1$. So, after sorting, we discover that
$z_{m_2}=w_{m_2};\,z_{m_2+1}=w_{m_2+1};\,\ldots,z_{k+1}=w_{k+1}$.
So in tableaux $T,T'$, the columns $m_1,\ldots,k$ agree as
does the bottom entry of column $k+1$. But the top entry in
column $k+1$ is larger in $T'$ than in $T$. Hence
$c_{T'}$ is lexically larger than $c_T$. 

\medbreak

Suppose now that $T'$ appears in $B$ in equation~\ref{case2:spanning:intoAB}.
Define $c_4=\min_{c_2\le i\le\lambda_1\atop w_{c_2}<v_i} i$.
Since $w_{c_2}\gp v_{c_2-1}=v_{c_0}$, we have
\[ T'=\rsktblx n{
\row{
v_{m_1},        , \ldots, v_{c_2-1}, \rlap{\kern-6pt$\srow{
w_{c_2}=,..  \ldots,. =w_{c_2}\makebox[0pt][r]{
\raisebox{10pt}{{Column:} $\sst c_4+s-1$}}
,..         \aword{\aword{\aword{W'}}},   } $}
}\\
  \smush{\row{
       ,        ,       ,    ,  ,  , ,
          \kern-25pt\raisebox{-4.75pt}{\makebox[0pt][l]{\frame{\row{   , , ,..... }}}}}}
  \row{
       , w_{m_2}, \ldots, w_{c_2-1}, \rlap{\kern-10pt$\srow{ 
w_{c_2+s}, \ldots, w_{c_3}, {\aword{W''}},  w_{c_3+1}, \ldots, w_{\lambda_2}}$}
         }
\row{ , , , , , ,.. }\\
  }, \]
where $W',W''$ is a shuffle of $v_{c_4}\cdots v_{\lambda_1}$ into
parts of size $\lambda_1-c_4-s+1$ and $1\le s\le c_3-c_2+1$ respectively
and, as before, the boxed elements must be sorted so that $T'$ will be row standard.

The top rows of $T,T'$ agree through
column $c_2-1$. Again either the bottom of column $c_2$ increases 
 and we are done or the
the number of positive letters in the bottom row that equal $w_{c_2}$ decreases.  
So iterating the straightening law on $T'$
 eventually increases the modified column word.

\end{proof}

The algorithm {\tt row-straighten} specializes to the
classical straightening for skew and partition shaped
tableau when $m_1\ge m_2$. The preceding result implies that
the column word also increases in the skew case.

\begin{corollary}\label{skew-strtning}
Proposition~\ref{lowerColWord} also holds with $m_1\ge m_2$.
\end{corollary}

\begin{proof}
We use what is often called the method of fake letters.
Fill cells $m_2-1,\ldots,m_1-1$ in the the top row with new negative letters disjoint 
from and smaller than the letters in $\C L$, we will name these letters 
$f_{m_2-1},\ldots,f_{m_1-1}$. 
Straighten this new tableau. 

The tableaux appearing in expression $B$ in the preceding proof
have the ``fake'' letters $f_{m_2-1},\ldots,f_{m_1-1}$ in the same positions as
does the original tableau $T$. 
If we apply the algebra homomorphism that sending $(f_i|j)$ to $\delta_{i,j}$,
$[T']$ is sent to $0$ for all $T'$ in the expression $A$ in the preceding proof
and the fake letters are erased from all other tableaux in the expression.
\end{proof}

We have now established the correctness of Algorithm~{\tt row-straighten}.
\begin{theorem}\label{basis}
The straight tableaux of shape $D$ form a $\Z$-basis for $\tabmod_D$
and Algorithm~{\tt straighten-tableau} expands any generator of 
$\tabmod_D$ in terms of this basis. Further, given a row-standard
tableau $T$, the expansion of $[T]$
is in terms of tableaux with larger column words than $w_T$.
\qed\end{theorem}

By Corollary~\ref{modincr} 
carefully analyzing the proof Proposition~\ref{lowerColWord} we can extend the preceding result.
\begin{corollary}\label{modincrtwo}
Algorithms {\tt row-straightening} and {\tt straighten-tableau}
produce tableaux with weakly larger modified column words than
that of the input tableau. \qed
\end{corollary}

\section{Flagged super-Schur modules}
\label{sect:flagged}

The {\em flagged} Schur modules $\tabmod_D_f$ have been the subject
of considerable interest
(see for instance \cite{LS90,RS96b,LM96}).
We apply the preceding results to flagged super-Schur modules of row-convex
shape. The fact that the initial terms in the straight basis are distinct allow
the straight bases to descend  to bases
of the corresponding flagged module.
I will start by formalizing the notion of a {flagged superSchur
module}. 
\begin{definition}
Let $\mbf f$ be a weakly increasing sequence of letters in the
alphabet $\C L$. Regard this sequence as indexed by elements
of $\C P$. The \defi{flagged superSchur
module} $S_{\mbf f}^D(\C L)$ is
the subquotient of $\Superlp LP$ equal to
the image of the
submodule $\tabmod_D(\C L)$ under the map $\phi_{\mbf f}$ which 
quotients $\Superlp LP$
by setting $(l|p)=0$ whenever $l>f_p$.

A tableau $T$ is \defi{flagged} if,  in each column~$i$,
$T$ has no entry exceeding $f_i$.
\end{definition}\index{{S_{\mbf f}^D(\C L)}}\index{{\phi_{\mbf f}}} 
If $T$ is row-standard (of any shape) and fails to be flagged, then
$\phi_f([T])=0$--
since  each monomial in the expansion of any row in which the
flagging condition is violated has some factor $(l|p)$ with $l>f_p$.

Classical results (see~\cite{Sta}) tell us that if $D$ is a skew-tableau then
a basis for $\tabmod_D_\mbf f$ is given by all $[T]$ such that $T$
is standard and flagged.
This result carries over to flagged row-convex superSchur modules.
\begin{theorem}
Let $D$ be a row-convex shape. Fix a weakly increasing flag $\mbf f$.
A basis for $\tabmod_D_\mbf f(\C L)$ is given by
the elements $\phi_\mbf f([T])$ where $T$ runs over all flagged, straight tableaux
of shape $D$ with entries chosen from $\C L$.
\end{theorem}
\begin{proof}
It suffices to show that the basis elements are linearly
independent, indeed that their initial terms under
any diagonal term order are still distinct.  This follows by  
Proposition~\ref{onesixsix}, and the  observation immediately
following it since the tableaux $T$ are both straight and flagged. 
\end{proof}

This result has the following easy generalization.  Let $\mbf f,\mbf g$
both be weakly increasing sequences of letters in $\C L$ indexed by elements
of $\C P$ such that $\mbf f\le \mbf g$ componentwise. Define the doubly flagged
superSchur module \indx{{\tabmod_D_{\mbf f,\mbf g}(\C L)}}
 to be the image of
$\tabmod_D(\C L)$ under the map \protect\indx{\phi_{\mbf f,\mbf g}} 
quotienting $\Superlp LP$ by the ideal generated by $\left\{ (l|p) :
l\not\in f_l,\ldots,g_l\right\}$. Call a tableau $T$ \defi{doubly flagged}
with respect to $\mbf f,\mbf g$ if every entry in column~$i$ is between
$f_i$ and $g_i$. The same proof as above shows the following.
\begin{theorem}
\label{doubflag}
A basis for $\tabmod_D_{\mbf f,\mbf g}(\C L)$ is given by
the elements $\phi_\mbf f([T])$ where $T$ runs over all 
doubly flagged, straight tableaux
of shape $D$ with entries chosen from $\C L$.
\qed\end{theorem}
When $\C L=\C L^-$ and $D$ is skew, this results appears in~\cite{Sta}.

\section{Groebner and {\protect\sc sagbi} bases.}
\label{sect:groebner}
The basis theorems developed above  have various ring-theoretic applications.
For convenience, we state them in terms of commutative rings,
i.e. we assume that the alphabet $\C L$ consists entirely of negative
letters. We can then take the subalgebra of the polynomial algebra 
$\k[x_{i,j}]$ generated by all polynomials indexed by tableaux of some
fixed shape.

\begin{definition}
Let $D$ be a shape.
Let $R^D$ be the subalgebra generated by all $[T]$ for all tableaux $T$ of
shape $D$. Allowing that $[T]$ has degree~1, this is a graded algebra generated
by its degree~1 part.
\end{definition}
These algebras turn out to be 
the homogeneous coordinate rings of certain {\em configuration varieties}
(see~\cite{M94} for example) embedded in projective space
by a combination of Plucker embeddings, Segre products and Veronese
maps. Configuration varieties parameterize a tuple of subspaces of $n$-space
subject to  certain lower bounds on the dimensions in which they may intersect.
The results in this section are stated for the rings $R^D$ but they hold equally
for the subquotient rings generated by the (doubly) flagged modules
$\tabmod_D_{\mbf f,\mbf g}(\C L)$.

Below we produce a Groebner basis for $R^D$. First we establish some notation.

\begin{definition}
Let $D$ be a shape. Let $T,T'$ be tableaux of shape $D$. Define $T\circ T'$ to
be the tableau (no longer of shape $D$) formed by alternating rows of $T$ and $T'$,
starting with the top row of $T$, then the top row of $T'$ etc.
\end{definition}

\noindent
Before stating the theorem we recall our convention that the initial
term of a polynomial is the {\em smallest} term in the polynomial and
define the degree of a Groebner basis be the 
highest degree of a polynomial appearing in Groebner basis. 

\begin{theorem}\label{GrBasis}
Let $D$ be an $n$-rowed, row-convex shape. 

Consider the polynomial ring whose variables consist of all straight tableaux of shape $D$ on $\C L$.
We define a graded term order on this ring as follows. We say $T'<T''$ when $c_{T'}<c_{T''}$ and we define
$\prod_{i=1}^k T'_i < \prod_{j=1}^k T''_i$ where $r<s$ implies $T'_r<T'_s$ and $T''_r<T''_s$
when $\mbf c_1(T'_k),\ldots, \mbf c_1(T'_1),\mbf c_2(T'_k),\ldots, \mbf c_2(T'_1),\ldots$ is smaller
in lexicographic order than $\mbf c_1(T''_k),\ldots, \mbf c_1(T''_1),\mbf c_2(T''_k),\ldots, \mbf c_2(T''_1),\ldots$
where $\mbf c_i(T)$ is the $i$th column of $T$ read from bottom to top.

If $I$ is the kernel of the ring map sending a straight tableau $T$ to $[T]$,
then $I$ has a degree~2 Groebner basis consisting of all
polynomials 
\begin{equation}\label{GroebnerRelns}
T'\cdot T'' \quad - \quad
\sum_\kappa  \beta_\kappa\cdot 
\begin{array}{c}
\\
\aword{{\SC t}'_1}\\
\aword{}\\
\aword{{\SC t}'_{i-1}}\\
\aword{v_\kappa}\\
\aword{{\SC t}'_{i+1}}\\
\aword{}\\
\aword{}\\
\aword{{\SC t}'_{n}}\\
\\
\end{array}
\,\cdot\,
\begin{array}{c}
\\
\aword{{\SC t}''_1}\\
\aword{}\\
\aword{}\\
\aword{{\SC t}''_{j-1}}\\
\aword{w_\kappa}\\
\aword{{\SC t}''_{j+1}}\\
\aword{}\\
\aword{{\SC t}''_{n}}\\
\\
\end{array}
\end{equation}
where $T',T''$ range over all straight tableaux of shape $D$ on $\C L$ and
where \[
\sum_\kappa \beta_\kappa\cdot 
    {\aword{\mbf{v_\kappa}} \atop \aword{\mbf{w_\kappa}}} =
     Syz_{a_1,\ldots,a_r;\, b_1,\ldots,b_s}\left(\aword{{\SC t}_i} \atop \aword {{\SC t}_j}\right)\,; \]
if $  {\aword{{\SC t}_i}\atop\aword{{\SC t}_j}} $ is not skew
$\mbf a$ and $\mbf b$ are allowed to range over all sequences of consecutive column 
indices such that $\mbf a$ ends in the last column of ${\SC t}'_{i}$ and
if $  {\aword{{\SC t}_i}\atop\aword{{\SC t}_j}} $ is skew, $\mbf a$ and $\mbf b$
range over all sequences of column indices such that $r+s$ is large than the number
of columns in $  {\aword{{\SC t}_i}\atop\aword{{\SC t}_j}} $.
\end{theorem}

\begin{proof}
Let $T\le T'$ be straight tableaux of shape $D$. 
Corollary~\ref{polzn-ident} verified that the polynomials~\ref{GroebnerRelns} are in $I$.
It suffices to show that if $T := T'\circ T''$ is
not straight then there exists a relation in the Groebner basis~\ref{GroebnerRelns} 
whose initial term is $T'\circ T''$.
Since $T$ is not straight, there exists $i\le j$ such that if $\SC t'_i$ is the $i$th
row of $T$ and $\SC t''_j$ is the $j$th row of $T'$ then the tableau 
$\aword{{\SC t}'_i}\atop \aword{{\SC t}''_j}$
is not straight. But then Proposition~\ref{lowerColWord} and Corollary~\ref{skew-strtning} show that
choosing $\mbf a$ and $\mbf b$ as in algorithm~{\tt row-straighten} gives $T'\cdot T''$ as the initial
term of the polynomial~\ref{GroebnerRelns}. 
\end{proof}

The Groebner basis of Theorem~\ref{GrBasis} is not reduced, nor are the initial terms of its 
elements polynomials necessarily distinct. The initial terms can be made distinct by choosing $\mbf a$
and $\mbf b$ as in algorithm {\tt row-straighten}. 
The above theorem does not require that $\C L=\C L^-$, although
one requires the notion of a non-commutative Groebner basis for general $\C L$.
Specifically when $\C L=\C L^-$, it is 
possible to restrict $r+s$ to be one more than the maximum of the number of columns in ${\SC t}_i$ and
the number of columns in ${\SC t}_j$ while simultaneously eliminating any other restrictions on
the indices $\mbf a$ and $\mbf b$; the underlying straightening law is presented in
Chapter~III of~\cite{thesis}. In general, that straightening law fails $\C L$ contains positive
letters.

\begin{porism}
The monomial $T_1\cdot T_2\ldots T_k$ is standard with respect to the Groebner basis in
Theorem~\ref{GrBasis} iff $T_1\circ T_2\circ \ldots\circ T_k$ is straight.
\qed\end{porism}

The existence of a degree~2 Groebner basis for an algebra is known to
imply that there is an (infinite) linear free resolution of the ground field over 
the algebra.

A {\sc sagbi} (Subalgebra Analogue of a Groebner Basis for Ideals)
basis, see~\cite{KM89} and~\cite{RoSw90}, is a
generating set for a subalgebra such that the
initial terms of the subalgebra are contained in
the algebra generated by the initial terms of the
generating set.

\begin{theorem}
Let $D$ be a row-convex shape.
The straight basis elements of shape $D$ form a {\sc sagbi} basis for 
$R^D$ with respect to any diagonal term order.
\end{theorem}
\begin{proof}
It suffices to show that if $p$ is in $R^D$ then its initial term is
the product of the initial terms of some multiset of straight tableaux.
Define $D^{\circ k}$ to be the shape formed by replacing each row of $D$ with $k$ copies
of itself.
So $\init_{}(p)$ is $\init_{}(p_k)$ where $p_k$ is the component of $p$ lying
in $S^{D^{\circ k}}$ with $k$ maximal such that $p_k\ne0$. Since the initial
terms of straight tableaux of fixed shape are distinct, $\init_{}(p)=[T]$
for some straight tableaux $T$ of shape $D^{\circ k}$. But we can write
$T=T_1\circ\cdots\circ T_k$. Each $T_i$ must be straight and
$\init_{}([T])=\prod_{i=1}^k \init_{}([T_i])$
\end{proof}

\begin{corollary}
Let $D$ be a row-convex shape.
The row-standard tableaux of shape $D$ form a {\sc sagbi} basis for 
$R^D$ with respect to any diagonal term order.
\qed\end{corollary}

By the usual results these {\sc sagbi} bases give algorithms for
determining whether a polynomial in variables $x_{i,j}$ belongs to
$R^D(\C L)$ or ({\em a forteriori}) $S^D(\C L)$ and, if so, writing
it in terms of the generators $[T]$ where $T$ is straight of shape $D$.
By results of~\cite{Stu96}, a {\sc sagbi} basis 
for an algebra allows that algebra to be deformed to an algebra generated
by monomials. In~\cite{thesis} and~\cite{straightApplic} this deformation is
used to prove that the subalgebra generated by all tableaux
of a fixed row-convex shape is Cohen-Macaulay.

\section{A branching rule and flagged corner-cell recurrence.}
\label{sect-br-rule}

Our final application concerns a branching rule for row-convex representations.
The Schur and Weyl modules $\tabmod_D(\C L^-)$ and $\tabmod_D(\C L^+)$ 
are $GL_n$ representations with $GL_n$-action induced by the algebra
homomorphism $g:\Superlp LP\rightarrow \Superlp LP$
given by $g\big((r|s)\big)=\sum_i g_{i,r} (i|s)$ where 
$g\in GL_n$ equals $(g_{i,j})$. In order to handle sets $\C L$ containing
letters of both positive and negative sign,  
we will work with representations of the general linear Lie superalgebra, $\gl_L$.

We express a $\gl_L$-representation, 
corresponding to a row-convex shape $D$, in terms of $\gl_{L\backslash \{a\}}$ 
representations (for some $a\in\C L$) corresponding to subshapes of $D$. 
The combinatorics for the case $\C L=\C L^+$ is identical to that of the
branching rule in~\cite{RS96c} and new when $\C L=\C L^-$. We present a filtration
that realizes this branching rule in a characteristic--free fashion. This provides the 
row-convex case of filtration conjectured in~\cite{RS96c} to exist for all \% comment deltd
It should be noted that the orientations of~\cite{RS96c} are at variance from those of~\cite{RS95};
we adhere to the orientation of the latter. Thus the term row-convex in~\cite{RS96c} should be read
as ``column-convex'' in the context of both~\cite{RS95} and the present paper. The branching rule 
presented below generalizes to the case of flagged super-Schur modules; branching rules for flagged
Schur modules are not treated in~\cite{RS96c}.

\addcontentsline{toc}{subsection}{Lie superalgebras.}

First we construct the general linear Lie superalgebras following Scheunert~\cite{Sc79}.

A free $\Z$-module $F$ is {\em signed}\glossary{signed module}
when it has distinguished
free submodules $F_0$ and $F_1$ whose direct sum is $F$. Elements of
$F_0$ and $F_1$ are called homogeneous and $|x|=i$ for $x\in F_i$.

A free signed $\integers$-module is a
\defi{Lie superalgebra} when it is endowed
with a \defi{superbracket} $[\, , \,]$ satisfying the commutativity relation,
\[
[x,y] = - (-1)^{|x||y|} [y,x]
\]
for homogeneous elements $x,y$ and the super-Jacobi identity
\[
(-1)^{|a||c|} [a,[b,c]] + (-1)^{|a||b|}[b,[c,a]] + (-1)^{|b||c|}[c,[a,b]]=0
\]
for homogeneous elements $a,b,c$.

Following~\cite{Brini-Teolis-J.Alg}, 
the general linear Lie superalgebra $\gl_L$, associated to the signed alphabet, $\C L$, 
\index{\gl_L} is the vector
space (over $\rationals$)  with basis $E_{a,b}$\index{E_{a,b}}
 for $a,b\in\cal L$, where $|E_{a,b}|=|a|+|b|$ and the 
bracket is
\[ [E_{a,b},E_{c,d}] = \delta_{b,c} E_{a,d} -
(-1)^{(|a|+|b|)(|c|+|d|)} \delta_{d,a} E_{c,b} . \] 
We next describe an action of $E_{a,b}$ on $\Superlp LP$.

\addcontentsline{toc}{subsection}{Polarizations}

A (left) \defi{superderivation} $D$ on a superalgebra $A$ is a
$\Z$-linear endomorphism of $A$ such that for $p,q$ homogeneous
in the $\Z_2$ grading of $A$, the identity
$D(pq)=(Dp)q+ (-1)^{\epsilon|p|}p(Dq)$ 
holds for some fixed $\epsilon\in\integers_2$. This $\epsilon$ is the 
{\em sign}\glossary{sign of a superderviation} of $D$, written $|D|$.

We define the \defi{letter polarization} 
$D_{a,b}:\SuperQlp LP \rightarrow \SuperQlp LP$ 
\index{D_{a,b}}
to be the superderivation with sign $|a|+|b|$
such that $D_{a,b} (c|p) = \delta_{b,c} (a|p)$ where
$\delta$ is the Kronecker delta.
It is easy to check that these superderivations are well-defined on
the $\Z$-subalgebra $\Superlp LP$.

The next example describes the action of the polarization operators 
in the case that the biproduct is the determinant of a minor.
\begin{example}\label{PolznNeg}
If $\C L=\C L^{-}$ and $\C P=\C P^{-}$, then $\Superlp LP$ is isomorphic to
$\Z[x_{i,j}\,:\, i\in\C L,\, j\in\C P]$. The action of $D_{i,j}$ (respectively $\R_{i,j}$
on this algebra is given by $\sum_{p\in\C P} x_{i,p}{\partial\over\partial x_{j,p}}$
(respectively $\sum_{l\in{\C L}} x_{l,i}{\partial\over \partial x_{l,i}}$).
\end{example}

To make our results characteristic free, we work over
$U(\gl_L)$,  the $\integers$-subalgebra of
the universal enveloping superalgebra of $\gl_L$
generated by all $E_{a,b}$, 
by $E_{a,b}^i\over i!$ for all $i\in\naturals$ and all
$a\ne b$ such that $|a|$=$|b|$, and by all  ${E_{a,a}\choose i}$, for
$i\in\naturals$.

\begin{proposition}
The map $E_{a,b}\mapsto D_{a,b}$, 
provides a representation 
of $U(\pl_L)$  on $\Superlp LP$.

If $D$ is a shape such that any letter
appearing in $\Deruyts(D)$ appears in $\C P$, then
this action descends to an action on $\tabmod_D(\C L)$.
\end{proposition}
{\em Proof Sketch. }

Defining the polarization operators to be superderivations on $\Super{\C L}$,
it can be shown that the polarization operators satisfy 
$D_{a,b}( p | q ) = (D_{a,b} p \,|\, q )$ where $(\,|\,)$ is the bilinear
form of Definition~\ref{bilinear}.
It is then clear that $\tabmod_D(\C L)$ is closed under the action of the
superderivations.

Details may be found in~\cite{thesis}.
\qed\smallskip

The branching rule for $\tabmod_D$ involves removing vertical or
horizontal strips from $D$.

\begin{definition}
 Let $D$ be a sorted row-convex shape. Define
a \defi{horizontal strip}, $E^+$, in $D$ to be any subset of
the cells of $D$ such that there exists a shape $D$ straight
tableau, $T$, on some alphabet $a^+<b_1<b_2<\cdots$
where the cells in $T$ that contain $a^+$ are precisely
the cells of $E$.
Similarly, define a \defi{vertical strip},  $E^-$ as any set of
cells containing all the negative letters $a^-$ appearing in some straight
tableau of shape $D$ on some alphabet $a^-<b_1<b_2<\cdots$.

Let $\mbf g, \mbf f$ be two weakly increasing sequences of letters
indexed by the elements of $\C P$. 
A vertical or horizontal
strip is $a$-flagged (with respect to $\mbf g,\mbf f$)
if it contains cells only in columns $i$
where $g_i\le a \le f_i$.
\end{definition}
Note that strips are allowed to be empty.

\begin{lemma}
Let $E$ be a vertical (respectively horizontal) strip in 
a row-convex shape $D$. Let ${\C I}_E$ be the multiset
of column indices indicating in which columns the cells of
the strip appear. If $E'$ is another vertical (respectively
horizontal) strip, then ${\C I}_E={\C I}_{E'}$ implies $E=E'$.
\end{lemma}
\begin{proof}
We utilize an alternative to Algorithm Straight-Filling for producing straight 
tableau with specified column content and shape. 
Suppose $D$ has $n$ cells. Define the desired contents of the columns of a 
tableau by a biword ${\mathbf u} = \left( \aword{\mbf{\hat u}} \atop \aword{\mbf{\check u}}\right)$
where $\hat\mbf u = w_{\Deruyts(D)}$ and 
$\check u_i \lm \check u_{i+1}$ if $\hat u_i = \hat u_{i+1}$.
Define a biword 
$\left( \aword{\mbf{\hat u}} \atop \aword{\mbf{\check u}}\right) =
\left( \hat u_{\sigma(1)},\ldots,\hat u_{\sigma(n)} \atop \check u_{\sigma(1)},\ldots,\check u_{\sigma(n)} \right)$
by permuting the entries of $\mathbf u$ so that $\check\mbf w$ weakly increases and so
$\check w_i=\check w_{i+1}$ implies $\hat w_i \le \hat w_{i+1}$; when $|\check w_i|=0$, this
inequality is strict.

Using this biword $\mathbf w$, we fill the tableau by starting with
an empty tableau of shape $D$ and adding successive letters reading left to right 
through the biword. At step~$j$ we place
$\check w_j$ in the northmost available cell (say row~$i$) in column~$\hat w_j$ such that
either $(i,\hat w_j-1)$ is not in the diagram or such that 
the cell $(i,\hat w_j-1)$ contains a letter $x$ with $x \lp \check w_j$. If
no such cell exists, then the biword does not arise from a straight tableau
of the given shape. To verify this algorithm, observe that if we
put $\check w_j$ into another row~$i'$, then either row-standardness is
violated or we have created a flippable inversion
in cells $(i,\hat w_j)$ and $(i',\hat w_j)$.

The lemma is an immediate consequence of the algorithm's correctness.
\end{proof}

\newcommand{\ch}{{ch}}
\begin{definition}
Suppose that $D$ is a row-convex shape, $\C L$ is
an alphabet, and $\Z[t_l\,:\, l\in\C L]$ is a polynomial ring.
Let $\mbf g, \mbf f$ be two weakly increasing sequences in $\C L$
indexed by the elements of $\C P$.  Define 
\[ \ch^D_{\mbf g,\mbf f}(\C L)=\sum_{T{\rm\ straight}} \prod_{(i,j)\in D} t_{T_{(i,j)}} \]
where the sum runs over all $\mbf g,\mbf f$-doubly flagged
straight tableaux of shape $D$ on $\C L$
and where $T_{(i,j)}$ is the $(i,j)$th entry of $T$. 
\end{definition}
When $\mbf f$ and $\mbf g$ are trivial, that is they contain respectively 
only the largest and smallest elements of $\C L$, and when $\C L$ contains 
letters of only one sign, this is the formal character of the
$GL({|\C L|})$-representation $\tabmod_D(\C L)$. If just one of $\mbf f,\mbf g$
is trivial, we get the formal character of a representation of a Borel subgroup.

The following identity is immediate from the definition of
a straight tableau, when $\C L=\C L^+$ it is due to~\cite{RS96c}.
\begin{proposition}
Fix two weakly increasing sequences $\mbf g,\mbf f$ of letters,
and choose $a\in\C L$.  If $D$ is a sorted row-convex diagram, then
\[
\ch^D_{\mbf g,\mbf f}(\C L)=
\sum_{E} \ch^{D/E}_{\mbf g,\mbf f}(\C L\backslash \{a\})
\]
where the sum runs over all $a$-flagged horizontal (respectively vertical)
strips $E$ in $D$ when $a\in\C L$ 
is positive (respectively negative), and where
$D/E$ is the diagram formed by removing $E$ from $D$.
\qed\end{proposition}

Preparatory to establishing a filtration for $\gl_L$-modules
$\tabmod_D(\C L)$ that realizes this identity
we define some components of that filtration. 
\begin{definition}
If $E, E'$ are two vertical (or two horizontal) strips in $D$,
define $E<E'$ in dominance order when for all $i$, the 
number (counted with multiplicity) of elements in $\{1,\ldots, i\}$
in $\C I_E$ is at least as large as the number of times these
elements appear in $\C I_{E'}$.

Let $E$ be a vertical (respectively horizontal) strip and let $a\in\C L$ be
negatively (respectively positively) signed. 
Define 
\[
\tabmod_{D,\ge E}(\C L;\, a)=\lspanZ_{T} \{[T]\}  
\]
where $T$ runs over all shape $D$ tableaux on $\C L$
in which $a$ appears in a vertical strip $E'$ 
weakly dominating $E$.

Define   $\tabmod_D{D, >E}(\C L;\, a)$
identically except for the requirement that $E'$ must strictly dominate $E$. 
\end{definition}
It is immediate from the definition that  $\tabmod_D{D, \ge E}(\C L;\, a)$
and  $\tabmod_D{D, >E}(\C L;\, a)$ are $\gl_{L\backslash\{a\}}$-representations.

\begin{theorem}
Let $D$ be a row-convex shape.
If $a$ is a negatively (respectively positively) signed letter in $\C L$
and $E$ is a horizontal (respectively vertical) strip in $D$, then
\[
\tabmod_{D,\ge E}(\C L;\, a)\biggm/ \tabmod_{D,> E}(\C L;\, a)
\quad\simeq\quad
\tabmod_{D/E} (\C L\backslash\{a\})
\]
as a $\gl_{L\backslash\{a\}}$-representation. Here $D/E$ is the shape formed by
removing $E$ from $D$.
\end{theorem}
\begin{proof}
It suffices to observe that given a row-standard
tableau $T$ such that the cells occupied by $a$ 
comprise $E$, then any tableaux appearing in the straightened
form of $[T]$ has the cells occupied by $a$ form a strip
$E'$ determined by a multiset $I'\ge I$. This can be seen by
directly examining the straightening relations. In particular,
any straightening relation which moves the $a$'s produces a
row-standard tableau in which the $a$'s form a horizontal
(respectively vertical) strip indexed by some $I'>I$.
\end{proof}

A more sophisticated result on the allowable contents of
a tableau appearing in the straightening of $[T]$ is proved
in~\cite{thesis} Chapter~III, Section~6.

\begin{corollary}\label{filtrn}
Let $D$ be a row-convex shape.
and let $a$ be a negatively (respectively positively) signed letter in $\C L$.
Let $E_1,\ldots,E_k$ be all vertical (respectively horizontal) strips in $D$
ordered compatibly with dominance, so that $i>j$ implies $E_j \not\ge E_k$.
The filtration
\[
\tabmod_D = 
\sum_{i=1}^k \tabmod_{D,E_i}(\C L;\, a) 
\supseteq\cdots\supseteq
\sum_{i=j}^{k} \tabmod_{D,E_i}(\C L;\, a) 
\supseteq\cdots\supseteq
\tabmod_{D,E_k}(\C L;\, a)  
\supseteq
0
\]
has $\tabmod_{D/E_j}(\C L\backslash\{a\})$
as the quotient, up to $\gl_{{\C L}\backslash\{a\}}$-isomorphism,
of its $j$th term by its $j+1$st term.
\qed\end{corollary}
If $\C L$ is sufficiently large, then the containments in the above filtration are all strict.

The preceding results generalize immediately to the $S_n$-representations  provided by
the Specht modules.

If we define $\tabmod_{D,E_i}_f(\C L;\, a)$ to be the flagged super-Schur
module found by taking the image of $\tabmod_{D,E_i}_f(\C L;\, a)$
under the map $(l|p)\mapsto 0$ when $l>f_p$, then we have the following.
\nonote{check proof}
\begin{proposition}
Maintaining the notation of Corollary~\ref{filtrn},
the filtration
\[
\tabmod_D_f = 
\sum_{i=1}^k \tabmod_{D,E_i}_f(\C L;\, a) 
\supseteq\cdots\supseteq
\sum_{i=j}^{k} \tabmod_{D,E_i}_f(\C L;\, a) 
\supseteq\cdots\supseteq
\tabmod_{D,E_k}_f(\C L;\, a)  
\supseteq
0
\]
of $\tabmod_D_f$ by $B$-modules
has $\tabmod_{D/E_j}_f(\C L\backslash\{a\})$
as the quotient, up to isomorphism,
of its $j$th term by its $j+1$st term.
Here $B$ is the subalgebra of $U(\pl_L)$
generated by all $E_{b,a}^i/i!$ for $b>a$ and
all $E_{a,a}\choose i$.
\qed\end{proposition}

When $\C L=\C L^-$, these results generalize to quantum Schur modules, details
appear in~\cite{thesis} and~\cite{quantStraight}.

\section{Acknowledgments}
I am indebted to David Buchsbaum for first asking whether the results of Woodcock
in~\cite{W94} could be strengthened to give a
canonical basis, to Mark Shimozono for pointing
out the connections to his work with Reiner in
general and to the corner-cell recurrence in
particular. This work first appeared as part of
the author's thesis~\cite{thesis}; I am deeply
grateful to my advisor Gian-Carlo Rota for his
long-standing support. This paper is dedicated to
his memory.

\done